\documentclass[11pt]{article}
\usepackage[margin=0.94in]{geometry}
\usepackage{epsfig,amssymb}
\usepackage[normalem]{ulem}
\usepackage{hyperref}
\usepackage{graphicx}
\usepackage{epstopdf}
\usepackage{amsmath}
\usepackage{bm}
\usepackage{setspace}
\usepackage{mathrsfs}
\usepackage{subcaption}
\usepackage[dvipsnames]{xcolor}
\usepackage{mathtools}
\usepackage{multirow}
\usepackage{lipsum}
\usepackage{comment}
\usepackage{enumerate}
\usepackage[affil-it]{authblk}
\usepackage{cite} 
\usepackage{leftidx}
\usepackage{relsize}
\usepackage{braket}
\usepackage{relsize}
\usepackage{float}
\usepackage{cite}
\usepackage{appendix}

\usepackage[autostyle]{csquotes}

\usepackage{fontenc}
\usepackage{comment}
\usepackage{soul}

\newcommand{\SAS}[1]{{\textcolor{black}{#1}}}
\newcommand{\etal}{\textit{et al}.~}
\providecommand{\keywords}[1]{\textbf{\textit{Keywords-}} #1}

\title{Thermodynamics of fractional-order nonlocal continua and its application to the thermoelastic response of beams}
\author[1]{Sai Sidhardh}
\affil[1]{School of Mechanical Engineering, Ray W. Herrick Laboratories, Purdue University, West Lafayette, IN 47907}
\author[1]{Sansit Patnaik}
\author[1]{Fabio Semperlotti$^\dagger$}

\begin{document}
\date{}
\maketitle

\begin{abstract}
This study presents a comprehensive framework for constitutive modeling of a frame-invariant fractional-order approach to nonlocal thermoelasticity in solids. For this purpose, thermodynamic and mechanical balance laws are derived for nonlocal solids modeled using the fractional-order continuum theory. This includes revisiting the Cauchy's hypothesis for surface traction vector in order to account for long-range interactions across the domain of nonlocal solid. Remarkably, it is shown that the fractional-order model allows the rigorous localized application of thermodynamic balance principles unlike existing integral approaches to nonlocal elasticity. Further, the mechanical governing equations of motion for the fractional-order solids obtained here are consistent with existing results from variational principles. These fractional-order governing equations involve self-adjoint operators and admit unique solutions, in contrast to analogous studies following the local Cauchy's hypothesis. To illustrate the efficacy of this framework, case-studies for the linear and the geometrically nonlinear responses of nonlocal beams subject to combined thermomechanical loads are considered here. Comparisons with existing integer-order integral nonlocal approaches highlight a consistent softening response of nonlocal structures predicted by the fractional-order framework, irrespective of the boundary and thermomechanical loading conditions. This latter aspect addresses an important incongruence often observed in strain-based integral approaches to nonlocal elasticity.\\

\noindent\keywords{Fractional calculus; Nonlocal elasticity; Constitutive Modeling; Geometric nonlinearity; Thermoelasticity}\\
$^\dagger$ All correspondence should be addressed to: \textit{fsemperl@purdue.edu}

\end{abstract}
\section{Introduction}
\label{sec: Introduction}

Several theoretical and experimental studies have shown that size-dependent effects, also referred to as nonlocal effects, are prominent in the response of complex structures of great relevance for many real-world applications. These size-dependent effects can be traced back to medium heterogeneity, existence of surface stresses, presence of thermal loads, and even medium geometry. More specifically, in the case of micro- and nano-structures, size-dependent effects have been traced back to the existence of surface and interface stresses due to nonlocal atomic interactions and Van der Waals forces \cite{sudak2003column,pradhan2009small,wang2011mechanisms}. In the case of macroscale structures, nonlocal effects can result from an ensemble of factors including material heterogeneity, interactions between layers (e.g. in FGMs or composite media) or unit cells (e.g. in periodic media), and geometric inhomogeneity \cite{romanoff2016using,hollkamp2019analysis,patnaik2019generalized}. In other terms, nonlocal governing equations for macrostructures often result from a process of homogenization of the initial inhomogeneous system. Further, geometric effects such as changes in curvature have also been shown to induce nonlocal size-dependent effects in nano-, micro-, and macro-structures \cite{sudak2003column,wang2011mechanisms}. 

The above mentioned complex slender structures have important applications in engineering and biotechnology. As an example, macroscale structures made from functionally graded materials (FGM) or sandwiched designs have been largely used in weight-critical applications such as aerospace, naval, and automotive systems \cite{kouchakzadeh2010panel,marzocca2011review}. Similarly, thin films, carbon nanotubes, monolayer graphene sheets and micro tubules have far-reaching applications in atomic devices, micro/nano-electromechanical devices, as well as sensors and biological implants \cite{sudak2003column,emam2013general}. Independently on the spatial scale, the key design constraints in the above applications include restrictions on space and weight. As a result, structural assemblies for lightweight applications are typically made of a combination of slender components like beams, plates, and shells. There are several applications where these structures are subject to large and rapidly varying mechanical and thermal loads that drive the system into a geometrically nonlinear regime.
A practical example includes the analysis of supersonic or hypersonic aerospace systems where the combination of large and quickly varying aero-thermomechanical loads induces highly a nonlinear response \cite{librescu2002supersonic,kouchakzadeh2010panel,marzocca2011review}. Similarly, the ability to account for coupled thermomechanical nonlinear effects is also critical in applications involving nano- and micro-structures such as, for example, in the design of biological implants, measurement devices, and sensors \cite{emam2013general,ebrahimi2015nonlocal}.
Despite the undeniable need for proper theoretical frameworks and computational tools capable of simulating the thermomechanical response of nonlinear and nonlocal structures, only a limited amount of studies focusing on geometrically nonlinear thermomechanical response of nonlocal slender structures are available in the literature. In the following, we briefly review the main characteristics of these studies and discuss key limitations.

Seminal works from Kro\"ner \cite{kroner1967elasticity} and Eringen \cite{eringen1972nonlocal} have explored the role of nonlocality in elasticity and laid its theoretical foundation. The key principle behind nonlocal theories relies on the idea that all the particles located within a prescribed area, typically indicated as the horizon of nonlocality, influence one another by means of long-range cohesive forces. This interaction between particles is accounted for by using gradient or integral relations for the strain field within the constitutive equations. These approaches lead to so-called “weak” gradient methods or “strong” integral methods, respectively. Integral methods \cite{eringen1972nonlocal,polizzotto2001nonlocal,barretta2018stress} capture nonlocal effects by re-defining the constitutive law in the form of a convolution integral of either the strain or the stress field over the horizon of nonlocality, whereas gradient elasticity theories \cite{peerlings2001critical,aifantis2003update,sidhardh2018inclusion,sidhardh2018exactisot} account for the nonlocal behavior by introducing strain or stress gradient dependent terms in the stress-strain constitutive law. As emphasized earlier, several applications involving nonlocal slender structures also experience thermal loads and geometric nonlinearities. Although several studies are available on the topics of geometrically nonlinear response of nonlocal slender structures \cite{yang2010nonlinear,srinivasa2013model,emam2013general} and on the thermomechanical response of nonlocal structures \cite{polizzotto2001nonlocal,shen2010nonlocal,tounsi2013nonlocal,ebrahimi2015nonlocal}, theoretical and numerical methods capable of addressing the combined geometrically nonlinear thermomechanical response of nonlocal structures have not been specifically addressed.

As mentioned previously, the classical studies on nonlocal (either linear or nonlinear) elasticity and nonlocal (linear) thermoelasticity encounter some key shortcomings. As an example, gradient theories experience difficulties when enforcing the boundary conditions associated with the strain gradient-dependent terms \cite{peerlings2001critical,aifantis2003update}. On the other side, the integral methods are better suited to deal with boundary conditions but they lead to mathematically ill-posed governing equations. This mathematical ill-posedness leads to erroneous predictions such as the absence of nonlocal effects or the occurrence of hardening behavior (not consistent with integral models) for certain combinations of boundary conditions \cite{romano2017constitutive}. In this class of problems, the ill-posedness stems from the fact that the constitutive relation between the bending moment field and the curvature is a Fredholm integral of the first kind, whose solution does not generally exists and, if it exists, it is not necessarily unique \cite{romano2017constitutive}. Additionally, in both these classes of methods, there are no available explicit relations to estimate the stress at a given point given the strain at that particular point. This latter aspect prevents the application of variational principles \cite{phadikar2010variational,anjomshoa2013application} and has critical implications on the development of a thermodynamic framework for the classical nonlocal approaches. More specifically, the modeling of nonlocality through nonlocal stress-strain constitutive relations allows only for a weak application (in a domain integral sense) and prevents the localized (point-wise) application of the thermodynamic balance laws. As discussed in \cite{polizzotto2001nonlocal}, the weak application of thermodynamic balance laws, particularly the second law, leads to inconsistencies in the nonlocal continuum framework.

In recent years, fractional calculus has emerged as a powerful mathematical tool to model a variety of nonlocal and multiscale phenomena. Fractional derivatives, which are a differ-integral class of operators, are intrinsically multiscale and provide a natural way to account for nonlocal effects \cite{podlubny1998fractional}. Given the multiscale nature of fractional operators, fractional calculus has found wide-spread applications in nonlocal elasticity \cite{cottone2009elastic,di2008long,carpinteri2014nonlocal,sumelka2014fractional,sumelka2015fractional,hollkamp2019analysis,patnaik2019generalized,hollkamp2020application}.
%
In a series of papers, Patnaik \etal \cite{patnaik2019generalized,patnaik2019FEM,sidhardh2020geometrically,patnaik2020plates,patnaik2020geometrically} have shown that a nonlocal continuum approach based on fractional-order kinematic relations provides an effective way to address the previously mentioned shortcomings of classical approaches to nonlocal elasticity. In addition, the formulation developed in these works is fully frame-invariant. Note that, unlike gradient elasticity methods, additional essential boundary conditions are not required when using the Caputo definition of the fractional operator \cite{hollkamp2019analysis,patnaik2019generalized}. Further, the nonlocal model based on fractional-order kinematic relations allows the application of variational principles and leads to well-posed governing equations that admit unique solutions \cite{patnaik2019FEM,patnaik2020plates}. \SAS{Finally we note that the fractional-order continuum theories are shown to be effective in developing reduced-order models for inhomogeneous systems, particularly periodic structures: periodic bar \cite{hollkamp2019analysis}, acoustic black hole structures \cite{hollkamp2020application}, and grid-stiffened plates \cite{patnaik2019generalized}.}

In this study, we build upon the fractional-order nonlocal continuum theory and develop a frame-invariant thermodynamic framework for the nonlocal solids. The overall goal of this study is two fold. 
First, to develop a thermodynamic framework for the fractional-order continuum formulation. We will show that the fractional-order continuum formulation allows for a rigorous application of all thermodynamic principles. More specifically, the use of fractional-order kinematic relations prevents the requirement of additional integral constitutive stress-strain relations as seen in classical nonlocal approaches (see, for example, \cite{polizzotto2001nonlocal}). {In fact, these additional residuals due to nonlocal interactions are assimilated within the corresponding fractional-order analogues. The most direct result is that the formulation does not require additional constraints associated with the thermodynamic balance laws.} It follows that the thermodynamic balance laws in the fractional-order theory are free from nonlocal residual terms, hence greatly simplifying the constitutive modeling of the nonlocal continuum and enabling a rigorous implementation of the thermodynamic principles at each point in the solid. The latter observation highlights an important benefit and a key motivation to pursue a fractional-order formulation to nonlocal thermoelasticity \SAS{as it sets the fractional approach apart from other classical nonlocal theories (either integral\cite{eringen1972linear,eringen1974theory} or differential\cite{eringen1983differential}) which typically can only satisfy thermodynamic equilibrium in a weak (global) sense.} Then, by enforcing the mechanical balance laws, we derive the governing equations of motion for the fractional-order nonlocal continuum. In this process, we employ a modification to Cauchy's hypothesis for the traction vector in order to include the effect of long-range interactions. \SAS{The thermodynamic framework is also of critical importance because it allows extending the energy-based methodologies typical of classical elasticity theories to the analysis of nonlocal structures. Examples of their applications include the linear buckling \cite{sidhardh2020fractional} and post-buckling \cite{sidhardh2020analysis} response of nonlocal structures that could not be tackled by standard integral models of nonlocal elasticity.}
The second objective of this study is to highlight the performance of the framework via a specific case-study focusing on a fractional-order Euler-Bernoulli beam subject to combined thermomechanical loads. For this purpose, we extend the fractional-order finite element model (f-FEM), originally developed in \cite{sidhardh2020geometrically}, to accurately solve the geometrically nonlinear fractional-order thermomechanical equations.

We note that several fractional-order thermomechanical models have been previously developed and presented in the literature \cite{povstenko2004fractional,povstenko2009thermoelasticity,povstenko2015fractional,vzecova2015heat}. However, these studies have focused primarily on the use of fractional-order operators to model complex thermal exchanges. More specifically, time-fractional operators have been employed to rewrite the heat conduction equation \cite{povstenko2004fractional,vzecova2015heat} in order to model dissipative effects associated with the thermal processes. Space-fractional operators have been used in \cite{povstenko2009thermoelasticity,vazquez2011fractional}, still within the heat transfer governing equation, to model anomalous forms of spatial diffusion. 
We emphasize that, differently from the present study, these previous works considered a local form of the stress-strain constitutive relation. On the contrary, our work does not consider fractional heat transfer, hence the heat conduction equation used in our study matches the classical integer-order form. In this regard, we merely note that a recent study \cite{vazquez2011fractional} has shown that the use of a space-fractional heat conduction equation leads to inconsistencies in the application of the second law of thermodynamics. This latter observation motivated us to use the classical (integer-order) heat conduction equation in our study.

The remainder of the paper is structured as follows: we begin with the development of a constitutive model for the fractional-order approach to nonlocal elasticity. This involves the thermodynamic and mechanical balance laws applied to the fractional-order solids. Later, to highlight the significance of the fractional-order models, we present the case-study of both a linear and a geometrically nonlinear response of nonlocal beams subject to combined thermomechanical loads and solved numerically via finite element techniques.

\section{Constitutive model for fractional-order thermoelasticity}
\label{sec: nonlocal_model}
In this section, we develop the constitutive model for fractional-order nonlocal thermoelasticity. As discussed earlier, the nonlocal beam theory presented in this study builds upon the formulation of a fractional-order nonlocal continuum presented in \cite{patnaik2019generalized}. \SAS{This formulation is a generalization of the seminal works on fractional-order continuum theories for nonlocal solids developed in \cite{drapaca2012fractional,carpinteri2014nonlocal,sumelka2014fractional}.} In the following, we review the key highlights of the continuum theory, and proceed with the development of constitutive model for nonlocal thermoelasticity.

Note the following notation used throughout the manuscript: $(\dot{\square})$ denotes the first integer-order derivative with respect to time, comma notation in the subscript $\square_{,j}$ will be used to denote integer-order spatial derivative with respect to the coordinate $x_j$, and Einstein summation is implied for repeated indices.

\subsection{Fundamentals of the fractional-order nonlocal continuum formulation}
\label{ssec: fundamental_frac}
Analogous to the classical approach to continuum mechanics, the response of a nonlocal solid can be analyzed by introducing two configurations, namely, the reference (undeformed) and the current (deformed) configurations. The motion of the body from the reference configuration (denoted as $\bm{X}$) to the current configuration (denoted as $\bm{x}$) is assumed as:
\begin{equation}
\label{eq: motion_description}
\bm{x}=\bm{\Phi}(\bm{X},t)
\end{equation}
such that $\bm{\Phi}(\bm{X},t)$ is a bijective mapping operation. The above mapping operation is used to model the differential line elements $\mathrm{d}\tilde{\bm{X}}$ and $\mathrm{d}\tilde{\bm{x}}$ in the undeformed and deformed configurations of the nonlocal solid using fractional-order operators. The fractional-order deformation gradient tensor $\overset{\alpha}{\textbf{F}}(\bm{x},\bm{X})$ defined with respect to nonlocal line elements is given by:
\begin{equation}
\label{eq: frac_def_tensor}
    \overset{\alpha}{\textbf{F}}=\frac{\mathrm{d}\tilde{\bm{x}}}{\mathrm{d}\tilde{\bm{X}}}
\end{equation}
In analogy with classical strain measures, the nonlocal strain can be defined using the fractional-order differential line elements as $\mathrm{d}\tilde{\bm{x}}\mathrm{d}\tilde{\bm{x}}-\mathrm{d}\tilde{\bm{X}}\mathrm{d}\tilde{\bm{X}}$. Following the above definition for fractional-order deformation-gradient tensor, the strain in the nonlocal solid is expressed as:
\begin{equation}
    \label{eq: frac_strain}
    \overset{\alpha}{\textbf{E}}=\frac{1}{2}\left(\overset{\alpha}{\textbf{F}}^T\overset{\alpha}{\textbf{F}}-\textbf{I}\right)
\end{equation}
Extending the above formalism, the Lagrangian strain tensor in the nonlocal medium is given by\cite{patnaik2019generalized,sidhardh2020geometrically}:
\begin{equation}
\label{eq: finite_fractional_strain}
\mathop{\textbf{E}}^{\alpha}=\frac{1}{2}\bigr(\underbrace{\nabla^\alpha {\textbf{U}}_X+\nabla^\alpha {\textbf{U}}_X^{T}}_{\tilde{\bm{\epsilon}}:~\text{Linear strain}}+\nabla^\alpha {\textbf{U}}_X^{T}\nabla^\alpha {\textbf{U}}_X\bigl)
\end{equation}
where $\textbf{U}(\bm{X})=\bm{x}(\bm{X})-\bm{X}$ denotes the displacement field. The fractional gradient denoted by $\nabla^\alpha\textbf{U}_X$ is given as $ \nabla^\alpha\textbf{U}_{X_{ij}} = D^{\alpha}_{X_j}U_i$ and consists of space-fractional derivatives. 
The space-fractional derivative $D^{\alpha}_{\bm{X}}\textbf{U}(\bm{X},t)$ is taken according to a Riesz-Caputo (RC) definition with order $\alpha\in(0,1)$ and it is defined on the interval $\bm{X} \in (\bm{X}_A,\bm{X}_B) \subseteq \mathbb{R}^3 $ in the following manner:
\begin{equation}
\label{eq: RC_definition}
	D^{\alpha}_{\bm{X}}\textbf{U}(\bm{X},t)=\frac{1}{2}\Gamma(2-\alpha)\big[\textbf{L}_{A}^{\alpha-1}~ {}^C_{\bm{X}_{A}}D^{\alpha}_{\bm{X}} \textbf{U}(\bm{X},t) - \textbf{L}_{B}^{\alpha-1}~ {}^C_{\bm{X}}D^{\alpha}_{\bm{X}_{B}}\textbf{U}(\bm{X},t)\big]
\end{equation}
where $\Gamma(\cdot)$ is the Gamma function, and ${}^C_{\bm{X}_{A}}D^{\alpha}_{\bm{X}}\textbf{U}$ and ${}^C_{\bm{X}}D^{\alpha}_{\bm{X}_{B}}\textbf{U}$ are the left- and right-handed Caputo derivatives of $\textbf{U}$, respectively. The complete expression for the nonlocal strain in Eq.~\eqref{eq: frac_strain} includes nonlinear terms required when accounting for large deformations. In this expression, we also highlight the linear component $\bm{\tilde{\epsilon}}$ of the fractional-order strain tensor that will be employed in studies of infinitesimal deformations in nonlocal solids. Unless otherwise specified, fractional-order strains in nonlocal solids will refer to the linear component $\bm{\tilde{\epsilon}}$. Before proceeding, it is worth discussing certain implications of this definition of the fractional-order derivative.
The interval of the fractional derivative $(\bm{X}_A,\bm{X}_B)$ defines the horizon of nonlocality (also called attenuation range in classical nonlocal elasticity). The length scale parameters $\textbf{L}_{A}^{\alpha-1}$ and $\textbf{L}_{B}^{\alpha-1}$ ensure the dimensional consistency of the deformation gradient tensor, and along with the term $\frac{1}{2}\Gamma(2-\alpha)$ ensure the frame invariance of the {strain-displacement} relations \cite{patnaik2019generalized,sumelka2014fractional}. \SAS{These length scales are independent parameters, which may or may not be equal to one another, hence resulting in an asymmetric horizon of nonlocality $(\bm{X}_A,\bm{X}_B)$ at the point $\textbf{X}$. The unequal values for the length scales ensure a truncation of the nonlocal region of influence for points close to the external boundaries and discontinuities within the solid. For complete details the reader should refer to \cite{patnaik2019generalized,patnaik2019FEM}. The choice of $\textbf{L}_A=\textbf{L}_B=l_f$ everywhere in the solid, where $l_f$ is a constant, identically reduces the above formulation to the fractional-order kinematic relations proposed in \cite{sumelka2014fractional}.}

In this formulation, nonlocality was introduced by using fractional-order kinematic relations. The fractional-order definition of the strain has critical implications on the thermodynamic framework for the fractional-order model of a nonlocal continuum. In the following section, we will show that the above approach to nonlocality allows the first and second law of thermodynamics to be enforced in a strong (localized) sense. In other terms, the fundamental laws of thermodynamics can be applied in a strict sense at each point in the nonlocal continuum as opposed to what happens in classical nonlocal approaches.

\subsection{Thermomechanical balance laws for fractional-order thermoelasticity}
\subsubsection{Thermodynamic balance laws}
\label{ssec: ThermoD_balance_laws}
In this section, we cast the fractional-order nonlocal model presented above within a thermodynamic framework. Consider a nonlocal solid $\Omega$ that undergoes the arbitrary motion $\bm{\Phi}$ which places a particle $p\in \Omega$ at $\bm{y}=\bm{\Phi}(p,t)$ at time $t$. At this instant, the solid occupies a domain $\mathcal{B}_t$ bounded by a surface $\partial \mathcal{B}_t$. Further, consider a part of the solid $\mathcal{P}\in \Omega$ that occupies the domain $\mathcal{D}_t=\bm{\Phi}(\mathcal{P},t)$ and such that point $p\in \mathcal{P}$. 

In order to enforce the first law of thermodynamics (i.e. the conservation of energy) over the arbitrary domain $\mathcal{D}_t$, we consider the energy balance applied to the domain $\mathcal{D}_t$. Later, we will show that the global balance laws hold true over any arbitrary domain $\mathcal{D}_t$ for the fractional-order continuum theory. This is unlike classical (integer-order) theories for nonlocal elasticity that allow the balance law to be applied only over the entire domain $\mathcal{B}_t$, that is in a global (weak) sense. The arbitrary choice of the domain, made possible by the fractional-order framework, also allows $\mathcal{D}_t$ to be reduced to an infinitesimal domain surrounding the point $p$ so to give the localized version of the thermodynamic law. The application of the first law of thermodynamics to the domain $\mathcal{D}_t$ allows writing the total energy $\mathcal{E}(\mathcal{D}_t,t)$ as the sum of the total internal energy $\mathcal{U}(\mathcal{D}_t,t)$ and the kinetic energy $\mathcal{T}(\mathcal{D}_t,t)$:
\begin{equation}
\label{eq: total_energy}
    \mathcal{E}(\mathcal{D}_t,t)=\underbrace{\int_{\mathcal{D}_t} \frac{1}{2} \rho(\bm{x},t)\dot{u}_{i}(\bm{x},t)\dot{u}_{i}(\bm{x},t)~\mathrm{d}V}_{\mathcal{T}(\mathcal{D}_t,t)}+\underbrace{\int_{\mathcal{D}_t} \rho(\bm{x},t)e(\bm{x},t)~\mathrm{d}V}_{\mathcal{U}(\mathcal{D}_t,t)}
\end{equation}
where $\rho$ is the mass density in the current configuration, $\dot{\textbf{u}}(=\dot{\bm{y}})$ is the spatial velocity of the particle at a generic point $p$ defined as the time derivative of the displacement field. In the above equation, $e(\bm{x},t)$ is the internal energy per unit mass, which is postulated to be a function of the mechanical and thermal state variables.

Recall that, for the classical (integer-order) local elastic solid, the internal energy of a point is a function of the integer-order strain and entropy defined at that point. For the nonlocal solid, in addition to the local strain energy, the internal energy at a point must also include the energy contributions from long-range cohesive forces by other points within the solid. Given the fractional-order kinematic relations described in \S\ref{ssec: fundamental_frac}, the contribution of the additional energy from nonlocal interactions is fully captured in the fractional-order strain $\tilde{\bm{\epsilon}}(\bm{x},t)$. 
It is immediate to see that the internal energy density evaluated at a point $\bm{x}$ is a function of the fractional-order strain ($\tilde{\bm{\epsilon}}(\bm{x},t)$). Consequently, we have the internal energy $e=e(\tilde{\bm{\epsilon}},\tilde{\eta})$ corresponding to the thermoelastic response defined entirely in terms of the state variables, that are the fractional-order strain $\tilde{\bm{\epsilon}}(\bm{x},t)$ and the entropy per unit mass $\tilde{\bm{\eta}}(\bm{x},t)$. This functional relationship is in sharp contrast with the thermodynamic framework for classical nonlocal approaches. 
More specifically, the thermodynamic framework for classical (integer-order) nonlocality leads to $e=e({\bm{\epsilon}},\mathcal{R}(\bm{\epsilon}),\tilde{\eta})$ where $\bm{\epsilon}(\bm{x},t)$ is the local strain field, and $\mathcal{R}(\bm{\epsilon})$ is a linear integral operator which models nonlocality in the solid \cite{eringen1974theory,balta1977theory,polizzotto2001nonlocal}. The additional functional relationship via $\mathcal{R}(\bm{\epsilon})$ is necessary in these approaches to account for energy contributions due to nonlocal interactions. To this regard, note that the fractional-order strain $\tilde{\bm{\epsilon}}$ combines the local integer-order strain $\bm{\epsilon}$ and its integral $\mathcal{R}(\bm{\epsilon})$ into a single term (see Eq.~(\ref{eq: RC_definition})) \cite{patnaik2019generalized}. Indeed, this is precisely the reason that allows expressing the internal energy density as $e=e(\tilde{\bm{\epsilon}},\tilde{\eta})$. It will be shown that the latter observation is significant as it allows the first law of thermodynamics to be applied in a strict sense at every point in the domain without additional constraints.
For this, we consider the first law of thermodynamics over the arbitrary domain $\mathcal{D}_t$:
\begin{equation}
    \label{eq: first_law_orig}
    \dot{\mathcal{E}}(\mathcal{D}_t,t)=\int_{\mathcal{D}_t} \rho(\bm{x},t)b_i(\bm{x},t)\dot{u}_{i}(\bm{x},t)~\mathrm{d}V+\int_{\partial \mathcal{D}_t} \tilde{\text{t}}_i(\bm{x},t,\bm{n})\dot{u}_i(\bm{x},t)~\mathrm{d}A+\int_{\mathcal{D}_t} \rho(\bm{x},t)r(\bm{x},t)~\mathrm{d}V-\int_{\partial \mathcal{D}_t} h(\bm{x},t,\bm{n})~\mathrm{d}A
\end{equation}
where the first two terms on the right side correspond to the mechanical power supplied to domain $\mathcal{D}_t$ by the body forces $\textbf{b}(\bm{x},t)$ and the surfaces forces applied per unit area $\tilde{\textbf{t}}(\bm{x},t,\bm{n})$, respectively. Here, $\bm{n}(\bm{x},t)$ denotes the outward normal to the surface $\partial\mathcal{D}_t$. Note that the surface forces on $\partial \mathcal{D}_t$ are applied by points external to the domain of interest. These surface forces include interactions within the infinitesimal neighborhood, as in (classical) local elasticity, and also the long-range interactions of the nonlocal solid. On the contrary, the body forces $\textbf{b}(\bm{x},t)$ applied at every point in the domain by external sources are local as they are independent of the nonlocal interactions within the solid. Additional terms in the above equation correspond to the rate of change in thermal energy of the body due to internal heat generation at the rate of $r(\bm{x},t)$ per unit volume, and the heat flux \textit{out} of the body at the rate of $h(\bm{x},t,\bm{n})$ through the surface $\partial \mathcal{D}_t$. Using Eq.~\eqref{eq: total_energy}, the first law of thermodynamics can be expressed as:
\begin{equation}
\label{eq: first_thermo_full}
\begin{split}
    \frac{\mathrm{d}}{\mathrm{d}t}\int_{\mathcal{D}_t} \rho(\bm{x},t) \left( \frac{1}{2} \dot{u}_{i}(\bm{x},t)\dot{u}_{i}(\bm{x},t)+e(\bm{x},t)\right)~\mathrm{d}V&=\int_{\mathcal{D}_t} \rho(\bm{x},t)b_i(\bm{x},t)\dot{u}_{i}(\bm{x},t)~\mathrm{d}V+\int_{\partial \mathcal{D}_t} \tilde{\text{t}}_i(\bm{x},t,\bm{n})\dot{u}_i(\bm{x},t)~\mathrm{d}A\\
    &+\int_{\mathcal{D}_t} \rho(\bm{x},t)r(\bm{x},t)~\mathrm{d}V-\int_{\partial \mathcal{D}_t} h(\bm{x},t,\bm{n})~\mathrm{d}A
\end{split}
\end{equation}
The above equation corresponding to the first thermodynamic balance law has important implications. Unlike analogous results from classical models of nonlocal elasticity, the above energy balance law holds true for arbitrary domain $\mathcal{D}_t$ within the nonlocal solid.
This integral form over the arbitrary domain $\mathcal{D}_t$ can be extended to develop the expressions for a localized imposition of the first law of thermodynamics at any point $p\in \mathcal{D}_t$, at any arbitrary time $t$. Recall that the internal energy density of a nonlocal solid must include contributions from long-range forces. For this purpose, modifications are introduced in classical integer-order models of nonlocal elasticity via the integral operator $\mathcal{R}(\mathbf{\epsilon})$ which restricts a localized imposition of the first thermodynamic law \cite{polizzotto2001nonlocal}. 

Before proceeding further, we briefly discuss key characteristics of the traction vector $\tilde{\textbf{t}}(\bm{x},t,\bm{n})$.
The surface traction at a point $\bm{x}$ is defined over an imaginary surface normal to the vector $\bm{n}(\bm{x},t)$ and passing through the point $\bm{x}$. This vector captures the forces acting on the surface due to the interaction of point $\bm{x}$ with the rest of the solid. Local elasticity which considers interactions restricted within an infinitesimal domain surrounding the point, employs the classical Cauchy's hypothesis to define the traction vector, that is $\textbf{t}=\bm{n}\cdot \bm{\sigma}$, where $\bm{\sigma}$ is the classical stress-tensor defined at the point of interest. In the fractional-order formulation, suitable modifications to the above definition are required to account for additional forces acting on the surface due to long-range interactions. For the fractional-order approach, we{ present }the following \SAS{generalization of the Cauchy's postulate}:\\

\noindent \textit{The surface traction $\tilde{\textbf{t}}(\bm{x},t,\bm{n})$ acting on an imaginary surface $\partial \mathcal{D}_t$ perpendicular to the normal vector $\bm{n}(\bm{x},t)$ and passing through the point $\bm{x}$ is given as:
\begin{equation}
    \label{eq: nonlocal_traxn}
    \tilde{\textbf{t}}(\bm{x},t,\bm{n})= {\bm{I}}_{\bm{n}}^{1-\alpha} \cdot \tilde{\bm{\sigma}} (\bm{x},t)
\end{equation}
}
where $\tilde{\bm{\sigma}}$ is the nonlocal stress evaluated at the point $\bm{x}$. The integral operator ${\bm{I}}_{\bm{n}}^{1-\alpha}$ defined over the horizon of influence of point $\bm{x}$ is given as:
\begin{equation}
\label{eq: integral_normal_operator}
    \bm{I}^{1-\alpha}_{\bm{n}}(\tilde{\bm{\sigma}}) = {I}^{1-\alpha}_{x_1}(n_1 \cdot\tilde{\bm{\sigma}}) \hat{e}_1 +  {I}^{1-\alpha}_{x_2}(n_2\cdot\tilde{\bm{\sigma}}) \hat{e}_2 +  {I}^{1-\alpha}_{x_3}(n_3\cdot\tilde{\bm{\sigma}}) \hat{e}_3
\end{equation}
such that $\bm{n}=n_i \hat{e}_i$ ($i=1,2,3$), $\hat{e}_i$ are the orthonormal basis vectors. Further, ${I}^{1-\alpha}_{i}(\cdot)$ is a Riesz-type fractional integral defined as: 
\begin{equation}
\label{eq: reisz integral_def}
    I^{1-\alpha}_{x_i} \chi =\frac{1}{2}\Gamma(2-\alpha) \left[ l_{B_{i}}^{\alpha-1} \left({}_{{x_i}-l_{B_{i}}}I_{{x_i}}^{1-\alpha} \chi \right) + l_{A_{i}}^{\alpha-1} \left({}_{{x_i}}I_{{x_i}+l_{A_{i}}}^{1-\alpha} \chi \right) \right]
\end{equation}
where ${}_{x_i-l_{B_{i}}}I_{{x_i}}^{1-\alpha}\chi$ and ${}_{{x_i}}I_{{x_i}+l_{A_{i}}}^{1-\alpha}\chi$ are the left and right fractional integrals to the order $\alpha\in (0,1)$ of an arbitrary function $\chi$, in the $\hat{e}_i$ direction. In this study, the definition of the normal vector is restricted to orthonormal triads of Cartesian coordinates. This is due to the limited developments in fractional vector calculus, as discussed in \cite{tarasov2008fractional,patnaik2020towards}. A detailed discussion on the above expression for the traction vector in a nonlocal solid is outlined in the Appendix.

Analogous to classical local elasticity, the traction vector defined in Eq.~\eqref{eq: nonlocal_traxn} is used within the global balance laws (Eq.~\eqref{eq: first_thermo_full}) to obtain the corresponding strong form. Additionally, we define the rate of heat transfer across the surface as $h(\bm{x},t,\bm{n})=\bm{q}(\bm{x},t)\cdot \bm{n}$, where the vector field $\bm{q}$ is referred to as the heat flux \textit{out} of the body. 
By substituting the definitions of the nonlocal traction and the heat flux into Eq.~\eqref{eq: first_thermo_full}, and then applying the divergence theorem gives:
\begin{equation}
    \label{eq: local_first_law_step1}
    \int_{\mathcal{D}_t} \left[ \rho \frac{De}{Dt}+\rho (\Ddot{{u}}_i-b_i)\dot{u}_i -\left(\bm{I}^{1-\alpha}_{n_i}\tilde{\sigma}_{ij}\dot{u}_j \right)_{,i} -\left(\rho r-q_{i,i} \right)\right] \mathrm{d}V=0
\end{equation}
Here, $D\cdot/Dt$ is the material derivative evaluated in the Eulerian configuration. By assuming continuous field variables and using the localization lemma we obtain:
\begin{equation}
    \rho \frac{De}{Dt}+\rho (\Ddot{{u}}_i-b_i)\dot{u}_i -\left(\bm{I}^{1-\alpha}_{n_i}\tilde{\sigma}_{ij}\dot{u}_j \right)_{,i} =\rho r-q_{i,i}~~\forall \bm{x}\in \mathcal{D}_t
\end{equation}
The above result clearly illustrates the possibility of enforcing the first law of thermodynamics (i.e. the energy balance) at any arbitrary point $p\in \mathcal{D}_t$ of the nonlocal solid. Note that, according to the fractional-order continuum theory, we do not require additional constraints to be imposed over the state variables. This is unlike the previous studies based on integer-order theories of nonlocal elasticity which required the energy balance to be complemented by an additional condition on the nonlocal residuals. Thus, a localized imposition of the energy balance law for nonlocal solids is derived here following the fractional-order models. In order to simplify the above expression to a more recognizable form, we use the mechanical balance laws for linear and angular momentum of the nonlocal solid, {derived} later in \S \ref{subsec: mech_bal_law} (refer Eqs. \eqref{eq: momentum_balance_loc} and \eqref{eq: ang_momentum_balance_loc}).
Employing the mechanical balance laws for the nonlocal solid, and following standard algebraic operations, we write:
\begin{equation}
    \label{eq: first_law_thermod}
    \rho\frac{De}{Dt}=\tilde{{\sigma}}_{ij}\dot{\tilde{\epsilon}}_{ij}+\rho r-q_{i,i}~~\forall \bm{x}\in\mathcal{B}_t
\end{equation}
Note that the strain tensor $\tilde{\epsilon}_{ij}$ in the above expression is evaluated using fractional-order derivatives defined in Eq.~\eqref{eq: frac_def_tensor}.
We emphasize again that the above localized form obtained for the fractional-order approach is in net contrast with classical nonlocal approaches where the conservation of the first law can only be applied in a weak sense (see, for example, \cite{polizzotto2001nonlocal}). 

Next, we apply the second law of thermodynamics to the fractional-order continuum model. Recall the Clausius-Duhem inequality applied over the entire solid $\mathcal{B}_t$ states:
\begin{equation}
    \label{eq: entropy_int_gen}
    \frac{d}{dt}\int_{\mathcal{B}_t}\rho(\bm{x},t)\tilde{\eta}_{\text{int}}(\bm{x},t)\mathrm{d}V=\frac{d}{dt}\int_{\mathcal{B}_t}\rho(\bm{x},t)\tilde{\eta}(\bm{x},t)\mathrm{d}V- \int_{\mathcal{B}_t}\frac{\rho(\bm{x},t)r(\bm{x},t)}{T(\bm{x},t)}+\int_{\partial \mathcal{B}_t}\frac{h(\bm{x},t,\bm{n}(\bm{x},t))}{T(\bm{x},t)}\mathrm{d}A\geq 0
\end{equation}
where $T$ denotes the temperature of the solid, and $\tilde{\eta}_{\text{int}}$ is the internal entropy production density. Recall that the localized form of the second law of thermodynamics states that the internal entropy production rate is non-negative for all points inside the solid, that is $D{\tilde{\eta}}_{\text{int}}/Dt\geq0~\forall~\bm{x}\in\mathcal{B}_t$. Classical approaches to nonlocal thermoelasticity satisfy this inequality only in a weak sense, that is the integral form given in the above equation \cite{eringen1972nonlocal,eringen1974theory,balta1977theory}. A detailed discussion of the consequent physical inconsistencies can be found in \cite{polizzotto2001nonlocal}. 

In analogy with the classical approach, we introduce the Legendre transformation $\psi=e-T\tilde{\eta}$, where $\psi$ denotes the Helmholtz free energy per unit mass. It follows that $\psi=\psi(\tilde{\bm{\epsilon}},T)$ which is different from classical nonlocal approaches wherein $\psi=\psi(\bm{\epsilon},\mathcal{R}(\bm{\epsilon}),T)$ \cite{polizzotto2001nonlocal}.
By using the Legendre transformation along with Eqs.~\eqref{eq: first_law_thermod} and \eqref{eq: entropy_int_gen}, we obtain the following local form of the second law of thermodynamics for the nonlocal solid:
\begin{equation}
    \label{eq: clausius_duhem}
    \rho T\frac{D{\tilde{\eta}}_{\text{int}}}{Dt}=\tilde{\sigma}_{ij}\dot{\tilde{\epsilon}}_{ij}-\rho \frac{D\psi}{Dt}-\rho \eta \dot{T}-T_{,i}\frac{q_{i}}{T}\geq 0
\end{equation}
Remarkably, the above fractional-order inequality matches, in its functional form, the classical Clausius-Duhem inequality. Equation~\eqref{eq: clausius_duhem} also highlights a clear difference compared with classical nonlocal formulations that require additional terms within the inequality as a result of the functional dependence of $\psi$ on $\mathcal{R}(\bm{\epsilon})$. As discussed in \cite{polizzotto2001nonlocal}, these additional terms within the inequality disappear only when a weak form is considered. However, as mentioned previously, satisfying the second law of thermodynamics only in a weak sense leads to nonphysical results. 

Thus, it appears that the fractional-order continuum theory for nonlocal solids allows the thermodynamic balance laws to be applied in a localized form; a key observation in order to establish the thermodynamic consistency of the fractional-order continuum theory and to derive rigorous constitutive models for the nonlocal solid. 
Note that the thermodynamic balance principles for the nonlocal solid have been presented in the Eulerian setting over the current configuration $\mathcal{B}_t$ of the solid at time instant $t$. However, assuming small deformation, the above results can be extended for the domain in reference configuration $\mathcal{B}_0$. For a generalized study of large deformations, resulting from both geometric and material nonlinearities simultaneously, this assumption will not be valid. In such cases, the Lagrangian analogues for the governing equations can be derived from the above results\cite{germain1983continuum}.

\subsubsection{Mechanical balance laws}
\label{subsec: mech_bal_law}
In the above discussion on thermodynamic balance laws, we have obtained the energy balance and entropy inequality laws for a fractional-order nonlocal solid. For the sake of completeness, in the following we also derive the mechanical balance laws corresponding to conservation of linear and angular momentum. For this purpose, we continue the discussion taking the point of view of the current configuration $\mathcal{B}_t$. Similar to the procedure outlined in \S\ref{ssec: ThermoD_balance_laws}, we begin with the integral balance laws over an arbitrary domain $\mathcal{D}_t$ and derive the localized forms of the governing equations for a point $p\in\mathcal{D}_t$. For the domain under consideration, the statement for the balance of linear momentum is:
\begin{subequations}
\begin{equation}
    \label{eq: momentum_balance_int}
    \frac{\mathrm{d}}{\mathrm{d}t}\int_{\mathcal{D}_t} \rho(\bm{x},t) \dot{\textbf{u}}(\bm{x},t)~\mathrm{d}V=\int_{\mathcal{D}_t}~\rho(\bm{x},t) \textbf{b}(\bm{x},t) \mathrm{d}V+\int_{\partial\mathcal{ D}_t}~\tilde{\textbf{t}}(\bm{x},t,\bm{n})~\mathrm{d}A
\end{equation}
and the balance of angular momentum in the absence of external couples is expressed as:
\begin{equation}
    \label{eq: ang_momentum_balance_int}
    \frac{\mathrm{d}}{\mathrm{d}t}\int_{\mathcal{D}_t} \bm{y}\times\rho(\bm{x},t) \dot{\textbf{u}}(\bm{x},t)~\mathrm{d}V=\int_{\mathcal{D}_t}~\bm{y}\times\rho(\bm{x},t) \textbf{b}(\bm{x},t) \mathrm{d}V+\int_{\partial\mathcal{ D}_t}~\bm{y}\times\tilde{\textbf{t}}(\bm{x},t,\bm{n})~\mathrm{d}A
\end{equation}
\end{subequations}
where the operator $'\times'$ denotes the exterior product. 
First, the above balance laws are simplified by imposing the classical (integer-order) result for the conservation of mass \cite{germain1983continuum}. 
Thereafter, by substituting the expression for surface traction $\tilde{\textbf{t}}$ given in Eq.~\eqref{eq: nonlocal_traxn} for the fractional-order solid, we obtain:
\begin{subequations}
\begin{equation}
    \label{eq: momentum_balance_int1}
    \int_{\mathcal{D}_t} \rho(\bm{x},t) \ddot{\textbf{u}}(\bm{x},t)~\mathrm{d}V=\int_{\mathcal{D}_t}~\rho(\bm{x},t) \textbf{b}(\bm{x},t) \mathrm{d}V+\int_{\partial\mathcal{ D}_t}~{\bm{I}}_{\bm{n}}^{1-\alpha} \cdot \tilde{\bm{\sigma}} (\bm{x},t)~\mathrm{d}A
\end{equation}
\begin{equation}
    \label{eq: ang_momentum_balance_int1}
     \int_{\mathcal{D}_t} \bm{y}\times\rho(\bm{x},t) \ddot{\textbf{u}}(\bm{x},t)~\mathrm{d}V=\int_{\mathcal{D}_t}~\bm{y}\times\rho(\bm{x},t) \textbf{b}(\bm{x},t) \mathrm{d}V+\int_{\partial\mathcal{ D}_t}~\bm{y}\times{\bm{I}}_{\bm{n}}^{1-\alpha} \cdot \tilde{\bm{\sigma}} (\bm{x},t)~\mathrm{d}A
\end{equation}
\end{subequations}
where we use: $\dot{\bm{y}}=\dot{\textbf{u}}$. Using the definition of the Reisz-integral within the definition of the surface traction and applying the divergence theorem we obtain \cite{patnaik2020towards}:
\begin{subequations}
\begin{equation}
    \label{eq: momentum_balance_int2}
    \int_{\mathcal{D}_t} \rho(\bm{x},t) \ddot{{u}}_j(\bm{x},t)~\mathrm{d}V=\int_{\mathcal{D}_t}~\rho(\bm{x},t) {b}_j(\bm{x},t) \mathrm{d}V+\int_{\partial\mathcal{ D}_t}~\left({{I}}_{n_i}^{1-\alpha} \tilde{{\sigma}}_{ij} (\bm{x},t)\right)_{,i}~\mathrm{d}V
\end{equation}
\begin{equation}
    \label{eq: ang_momentum_balance_int2}
     \int_{\mathcal{D}_t} \varepsilon_{ijk}{y}_{j}~\rho(\bm{x},t) \ddot{{u}}_{k}(\bm{x},t)~\mathrm{d}V=\int_{\mathcal{D}_t}~\varepsilon_{ijk}~{y}_{j}\rho(\bm{x},t) {b}_{k}(\bm{x},t) \mathrm{d}V+\int_{\partial\mathcal{ D}_t}~\varepsilon_{ijk}~\left({y}_{j}{\bm{I}}_{n_{m}}^{1-\alpha} \tilde{{\sigma}}_{mk} (\bm{x},t)\right)_{,m}~\mathrm{d}V
\end{equation}
\end{subequations}
Applying the localization lemma to Eq.~\eqref{eq: momentum_balance_int2}, we arrive at the following result:
\begin{equation}
    \label{eq: momentum_balance_loc}
    \rho(\bm{x},t) \ddot{{u}}_j(\bm{x},t)=\rho(\bm{x},t) {b}_j(\bm{x},t) \mathrm{d}V+ \mathfrak{D}_{x_i}^{\alpha}\tilde{{\sigma}}_{ij} (\bm{x},t),~~~~\forall \bm{x}\in \mathcal{B}_t
\end{equation}
which is the localized linear momentum balance law for fractional-order solids. Here, $\mathfrak{D}_{x_i}^{\alpha}(\cdot)$ is the Riesz Riemann-Liouville derivative of order $\alpha$ which is defined as \cite{patnaik2019FEM}:
\begin{equation}
    \label{eq: r_rl_frac_der_def}
    \mathfrak{D}^{\alpha}_{x_i}\chi = \frac{1}{2}\Gamma(2-\alpha) \left[ l_{B}^{\alpha-1} \left({}^{~~~RL}_{x_i - l_{B}} D^{\alpha}_{x_i} \chi\right) - l_{A}^{\alpha-1} \left( {}^{RL}_{x_i}D^{\alpha}_{x_i + l_{A}} \chi\right)\right]
\end{equation}
where $\chi$ is an arbitrary function and ${}^{~~~RL}_{x_i-l_{B}}D^{\alpha}_{x_i}\chi$ and ${}^{RL}_{x_i}D^{\alpha}_{x_i+l_{A}}\chi$ are the left- and right-handed Riemann Liouville derivatives of $\chi$ to the order $\alpha$, respectively. 
In the derivation of the localized form of the linear momentum balance law, given in Eq.~\eqref{eq: momentum_balance_loc} from Eq.~\eqref{eq: momentum_balance_int2}, we used the following relation:
\begin{equation}
    \mathfrak{D}^\alpha_{x_i} \chi = \frac{\mathrm{d}}{\mathrm{d}x_i} \left[ I^{1-\alpha}_{i} \chi \right]
\end{equation}
The governing equations of motion derived here using the balance principles agree very well with the results obtained via variational principles \cite{patnaik2020towards}. 
Using the above result for linear momentum balance and $y_{i,j}=\delta_{ij}$ ($\delta_{ij}$ is the Kr\"onecker delta), the strong form of the angular momentum balance for fractional-order solids reduces to the symmetry condition of the nonlocal stress tensor:
\begin{equation}
\label{eq: ang_momentum_balance_loc}
    \tilde{\sigma}_{ij}=\tilde{\sigma}_{ji},~~~~\forall \bm{x}\in \mathcal{B}_t
\end{equation}

It is interesting to note the fractional-order Riesz-type Riemann-Liouville fractional-order derivative (divergence) of the nonlocal stress tensor in the elastodynamic equation (Eq.~\ref{eq: momentum_balance_loc}). This result differ from existing studies on fractional-order continuum theories that employ a first-order derivative of the nonlocal stress ($D^{1}_{x_i}\tilde{\sigma}_{ij}$) tensor in the same equation (see, for example, \cite{sumelka2014fractional}). This difference arises due to the consideration of long-range interactions within the expression for traction given in Eq.~\eqref{eq: nonlocal_traxn}. The additional contributions to the surface traction at a point in the nonlocal solid follow from the long-range interactions (see Appendix). Note that the Riesz-type Riemann-Liouville operator in the mechanical governing equations given in Eq.~\eqref{eq: momentum_balance_loc} is self-adjoint and the system is positive-definite. This result was established and proved in \cite{patnaik2019FEM}. These observations are clearly in contrast with either classical integer-order approaches to nonlocal elasticity, that have shown it is not possible to define a self-adjoint quadratic potential energy \cite{reddy2010nonlocal,challamel2014nonconservativeness}, or some fractional-order models  \cite{sumelka2014fractional,sumelka2015fractional}. Given the self-adjoint and positive-definite nature of our formulation, the resulting system of equations is well-posed and admit a unique solution \cite{patnaik2019FEM}. We will show in \S\ref{sec: num_results} that this well-posedness results in a consistent softening behavior of the structure {with increasing degree of nonlocality} irrespective of the thermomechanical load distributions and boundary conditions. {This result is significant because it bypasses a key inconsistency observed in classical nonlocal models and associated with the non self-adjointness of the operators \cite{reddy2010nonlocal,challamel2014nonconservativeness}.}

\subsection{Constitutive framework for fractional-order thermoelasticity}
\subsubsection{Constitutive modeling}
\label{subsec: constt_model}
The inequality in Eq.~(\ref{eq: clausius_duhem}) is used to derive the thermodynamically-consistent constitutive equations for fractional-order nonlocal elasticity. By substituting the expression for the time derivative of the Helmholtz free energy, the inequality in Eq.~(\ref{eq: clausius_duhem}) is expressed as:
\begin{equation}
    \label{eq: constt_1}
   T\dot{\tilde{\eta}}_{\text{int}}= \left(\tilde{\sigma}_{ij}-\rho\frac{\partial \psi}{\partial \tilde{\epsilon}_{ij}}\right)\dot{\tilde{\epsilon}}_{ij}-\rho\left(\eta+\frac{\partial \psi}{\partial T}\right)\dot{T}-T_{,i}\frac{q_{i}}{T}\geq 0
\end{equation}
Since the above inequality must hold for all thermoelastic processes as well as for arbitrary choices of the independent fields $\dot{\tilde{\epsilon}}_{ij}$ and $\dot{{T}}$, we obtain the following constitutive laws:
\begin{subequations}
    \label{constt_eq}
\begin{equation}
\label{eq: mech_constt_eq}
    \tilde{\sigma}_{ij}={\rho}\frac{\partial \psi}{\partial \dot{\tilde{\epsilon}}_{ij}}, ~~\forall~ \bm{x} \in \mathcal{B}_t
\end{equation}
\begin{equation}
\label{eq: therm_constt_eq}
    \tilde{\eta}=-\frac{\partial \psi}{\partial T}, ~~\forall ~\bm{x} \in \mathcal{B}_t
\end{equation}
\end{subequations}
Under assumptions of linear elasticity, the above equations are the fractional analogues for the Duhamel-Neumann's laws for classical thermoelasticity. Further, by using the above constitutive relations within Eq.~(\ref{eq: clausius_duhem}), the inequality reduces to:
\begin{equation}
\label{eq: reduced_ineq}
    \rho T\dot{\tilde{\eta}}_{\text{int}}=-T_{,i}\frac{q_{i}}{T}\geq 0, ~~\forall ~ \bm{x} \in \mathcal{B}_t
\end{equation}
which establishes the second law of thermodynamics for a fractional order nonlocal solid. The relations in Eqs.~(\ref{constt_eq},\ref{eq: reduced_ineq}) can be expressed as:\newline

\noindent\textbf{{Theorem}}: \textit{The constitutive relations for fractional-order nonlocal thermoelasticity do not violate the Clausius-Duhem inequality if they are of the form given in Eq.~(\ref{constt_eq}) and subject to Eq.~\eqref{eq: reduced_ineq}.}\newline

A few additional comments on this thermodynamic framework are needed. At a first glance, the form of the stress-strain constitutive relation in Eq.~(\ref{eq: mech_constt_eq}) might be deceiving as it appears to lead to a classical constitutive relation. 
Although this is, formally, a correct statement it does not entirely capture the nature of Eq.~(\ref{eq: mech_constt_eq}). As highlighted earlier, nonlocality was modeled using fractional-order kinematic relations given in Eq.~(\ref{eq: finite_fractional_strain}). Therefore, the stress defined through the Eq.~(\ref{eq: mech_constt_eq}) is also nonlocal in nature. 
In addition, this construction of nonlocality (i.e. based on fractional-order kinematic relations) allows the application of variational principles, ensures a quadratic form of the potential energy of the system, and leads to well-posed nonlocal governing equations \cite{patnaik2019FEM,patnaik2020plates}. 

Note also that, in the above study, we assumed an integer-order Fourier heat conduction law, that is $q_i=-kT_{,i}$, where $k$ is the material conductivity constant such that $k>0$. It is immediate that the inequality in Eq.~(\ref{eq: reduced_ineq}) obtained from the second law is trivially satisfied for the integer-order heat conduction law. 
Finally, following from the latter remark, a space-fractional thermal conduction law defined as $q_i=-k{D}_{X_i}^{\alpha}\theta$ was proposed in \cite{povstenko2009thermoelasticity}. However, as shown in \cite{vazquez2011fractional}, the space-fractional heat conduction law violates the second law of thermodynamics. Thus, we limit the scope of the current study to fractional-order constitutive modeling for the mechanical field but integer-order models for the thermal fields. 

\subsubsection{Linear fractional-order thermoelasticity}

In this section, we derive the constitutive relations for linear thermoelastic response of fractional-order nonlocal solids. For this purpose, recalling that the fractional-order nonlocal formulation allows a localized implementation of the thermodynamic principles, we write the Helmholtz free energy density for the thermoelastic response following the typical approach for local elasticity (albeit by using the fractional-order strain). Further, we use the Helmholtz free energy density to construct the material constitutive relations for nonlocal thermoelasticity. 
In this study on linear elastic behavior of nonlocal solids using fractional-order theories of thermoelasticity,  we extend the linear material constitutive relations to the geometrically nonlinear response of the fractional-order nonlocal solids (Saint Venant-Kirchhoff material model). Therefore, assuming small deformations, the distinction between reference and current configurations for the domain $\Omega$ vanishes. This implies the constitutive relations given in Eq.~\eqref{constt_eq} for domain $\Omega$ hold true in all configurations ($\mathcal{B}_0$ \& $\mathcal{B}_t$). 

We can cast the constitutive relations in Eq.~\eqref{constt_eq} in the form:
\begin{equation}
    \tilde{\sigma}_{ij}=\frac{\partial \mathcal{W}}{\partial \tilde{\epsilon}_{ij}},~~~~~~\tilde{\eta}=-\rho_0^{-1}\frac{\partial \mathcal{W}}{\partial T},~~~\forall \bm{X}\in \Omega
\end{equation}
where the total free energy $\mathcal{W}=\rho_0\psi$ is expressed in terms of the mass density $\rho_0$ in reference configuration. We write the free energy for an isotropic material as a series expansion of the fractional-order strain $\tilde{\epsilon}_{ij}$ and the temperature difference $\theta=T-T_0$, which is the difference between the temperature $T$ at any point within the continuum and the uniform ambient temperature $T_0$ at the reference state. The free energy, expanded in power series with respect to the strain invariants and the temperature difference, is given as \cite{oden1969finite}:
\begin{subequations}
\begin{equation}
    \label{eq: helm_general}
    \begin{split}
\mathcal{W}=a_0+a_1\tilde{J}_1+a_2\tilde{J}_2+a_3\tilde{J}_3+a_4\theta+a_5\tilde{J}_1^2+a_6\tilde{J}_2^2+a_7\tilde{J}_3^2+a_8\tilde{J}_1\tilde{J}_2&\\
    +a_9\tilde{J}_1\tilde{J}_3+a_{10}\tilde{J}_2\tilde{J}_3+a_{11}\theta^2+a_{12}\tilde{J}_1\theta+a_{13}\tilde{J}_2\theta+a_{14}\tilde{J}_3\theta+\text{h.o.t}&
    \end{split}
\end{equation}
where $a_k$ are material constants and
\begin{equation}
    \tilde{J}_1=\tilde{\epsilon}_{ii}; ~~~\tilde{J}_2=\frac{1}{2}\left(\tilde{\epsilon}_{ii}\tilde{\epsilon}_{jj}-\tilde{\epsilon}_{ij}\tilde{\epsilon}_{ij}\right);~~~\tilde{J}_3=\text{det}(\tilde{\epsilon}_{ij})
\end{equation}
\end{subequations}
are the invariants of the nonlocal strain tensor $\tilde{\bm{\epsilon}}$. Assuming that the solid is stress free in the undeformed state and the free energy is restricted to linear isotropic thermoelasticity (i.e. ignoring the higher order terms in Eq.~(\ref{eq: helm_general})), we obtain the following expression for $\mathcal{W}$:
\begin{equation}
\mathcal{W}=a_2\tilde{J}_2+a_5\tilde{J}_1^2+a_{11}\theta^2+a_{12}\tilde{J}_1\theta
\end{equation}
where the material constants are given as \cite{kelly2020material}:
\begin{equation}
    a_2=-2\mu;~~~a_5=\frac{1}{2}\left(\lambda+2\mu\right);~~~a_{11}=-\frac{\rho_0 C_v^0}{2T_0};~~~a_{12}=(3\lambda+2\mu)\alpha_0
\end{equation}
The material constants $\lambda$ and $\mu$ are the isothermal Lam\'e parameters for the isotropic solid, $\alpha_0$ is the coefficient of volumetric thermal expansion and $C_v^0$ is the specific heat at constant strain. The parameters $\alpha_0$ and $C_v^0$ are all defined in the reference state at $T_0$. Thus, the Helmholtz free energy density (per unit volume) for the thermoelastic response of a nonlocal isotropic solid is given by \cite{washizu1975variational,kelly2020material}:
\begin{equation}
    \label{eq: helmholtz_energy}
    \mathcal{W}=\frac{1}{2}\lambda\tilde{\epsilon}_{kk}\tilde{\epsilon}_{ll}+\mu\tilde{\epsilon}_{ij}\tilde{\epsilon}_{ij}-(3\lambda+2\mu)\alpha_0\tilde{\epsilon}_{kk}\theta-\frac{\rho_0 C_v^0}{2T_0}\theta^2
\end{equation}
By using the above expression for $\psi$ together with Eq.~\eqref{constt_eq}, the thermoelastic constitutive relations relating the different physical quantities for the isotropic solid are obtained as:
\begin{subequations}
\label{eq: constt_gen}
\begin{equation}
    \label{eq: constt_isot}
    \tilde{\sigma}_{ij}=\lambda\delta_{ij}\tilde{\epsilon}_{kk}+2\mu\tilde{\epsilon}_{ij}-(3\lambda+2\mu)\alpha_0\delta_{ij} \theta 
\end{equation}
\begin{equation}
    \tilde{\eta}=\rho_0^{-1}(3\lambda+2\mu)\alpha_0\tilde{\epsilon}_{kk}+\frac{C_v^0}{T_0}\theta
\end{equation}
\end{subequations}
Using the above results, the Helmholtz free energy in Eq.~\eqref{eq: helmholtz_energy} is recast in the following manner:
\begin{equation}
    \label{eq: helmoltz_simp}
    \mathcal{W}=\frac{1}{2}~\tilde{\sigma}_{ij}~\tilde{\epsilon}_{ij}-\frac{1}{2}~\rho_0~\tilde{\eta}~\theta
\end{equation}

\section{Thermoelastic Euler-Bernoulli nonlocal beam model}
In this section, we use the thermoelastic constitutive relations developed for the nonlocal solid to analyze the thermoelastic response of a fractional-order Euler-Bernoulli beam. 
Building on \cite{patnaik2019FEM,sidhardh2020geometrically}, we derive the geometrically nonlinear governing equations and the corresponding boundary conditions for the thermoelastic boundary value problem (BVP) using variational principles.

\subsection{Geometrically nonlinear constitutive relations}
\label{sec: nonlinear_constt}
\begin{figure}[h]
    \centering
    \includegraphics[width=0.6\textwidth]{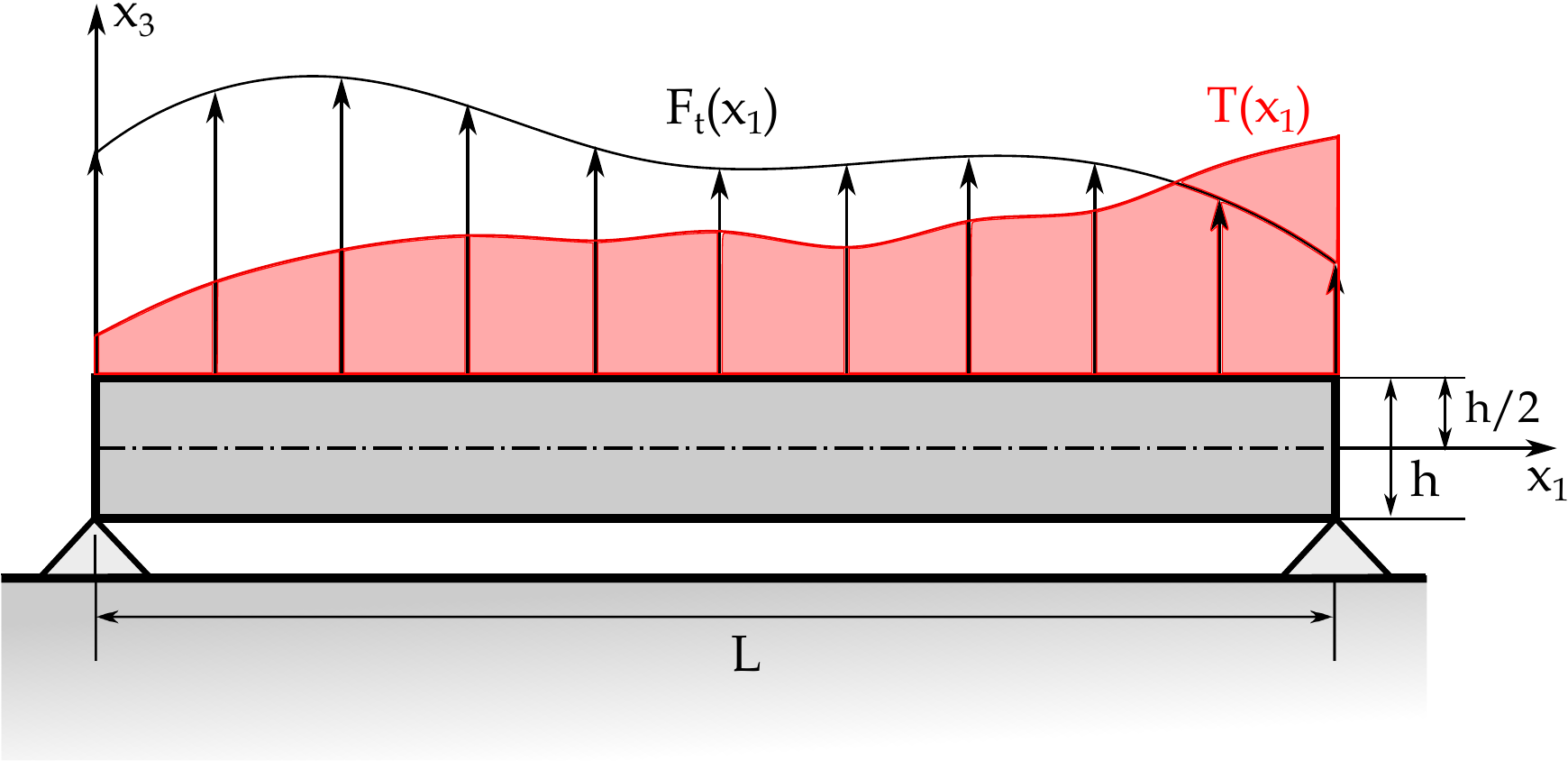}
    \caption{Schematic of an elastic beam subject to distributed transverse mechanical load $F_t(x_1)$ and thermal load $T(x_1)$.}
    \label{fig: schematic_beam}
\end{figure}

Consider a nonlocal beam subject to distributed transverse mechanical and thermal loads as illustrated in Fig.~(\ref{fig: schematic_beam}). As indicated in the schematic, the Cartesian coordinates for the current study are chosen such that $x_3=\pm h/2$ coincides with the top and bottom surfaces of the beam, and $x_1=0$ and $x_1=L$ are the ends of the beam along the longitudinal direction. The surface $x_3=0$ coincides with the mid-plane of the beam and the origin of the reference frame is chosen at the intersection of the mid-plane with the left-end of the beam.

The axial and transverse components of the displacement field $\textbf{u}(x_1,x_3)$ are denoted by $u_1(x_1,x_3)$ and $u_3(x_1,x_3)$, respectively. These displacement fields are given by the Euler-Bernoulli theory as:
\begin{subequations}
    \label{eq: euler_bern_beam}
    \begin{equation}
        u_1(x_1,x_3)=u_0(x_1)-x_3\left[\frac{\mathrm{d}w_0(x_1)}{\mathrm{d}x_1}\right]
    \end{equation}
    \begin{equation}
        u_3(x_1,x_3)=w_0(x_1)
    \end{equation}
\end{subequations}
where $u_0(x_1)$ and $w_0(x_1)$ are the mid-plane axial and transverse displacements, respectively. For a geometrically nonlinear analysis, assuming moderate rotations ($10^\circ-15^\circ$) but small strains, the fractional-order Lagrangian strain tensor in Eq.~\eqref{eq: finite_fractional_strain} can be further simplified using von-K\'arm\'an relations. The resulting fractional-order von-K\'arm\'an strain-displacement relations are given as \cite{sidhardh2020geometrically}:
\begin{equation}
    \label{eq: vonkarman_strain}
    \tilde{\epsilon}_{11}(x_1,x_3) = D_{x_1}^\alpha u_1(x_1) + \frac{1}{2}\left[D_{x_1}^\alpha u_3(x_1)\right]^2
\end{equation}
where $D^\alpha_{x_1}(\cdot)$ denotes the fractional-order RC derivative defined in Eq.~(\ref{eq: RC_definition}). {Note that $\tilde{\epsilon}_{11}$ here is nonlinear in transverse displacements. }As discussed previously, for the RC derivative used above, the 1D domain $(\bm{x}_A,\bm{x}_B)$ along the mid-plane of the beam is the horizon of nonlocal interaction at $\bm{x} (x_1,0)$. The end-points of the nonlocal horizon $\bm{x}_A (x_{A_1},0)$ and $\bm{x}_B(x_{B_1},0)$ are the terminals of the left- and right-handed Caputo derivatives within the RC derivative. It follows from Eq.~(\ref{eq: RC_definition}) that $l_A=x_1-x_{A_1}$ and $l_B=x_{B_1}-x_1$ are the length scales along $\hat{\bm{x}}_1$ to the left and right hand sides of the point $\bm{x}(x_1,0)$, respectively. \SAS{For the Euler-Bernoulli beam displacement field given in Eq.~\eqref{eq: euler_bern_beam}, a non-zero expression for the transverse shear strain would be obtained. However, for the slender beam assumed here the rigidity to transverse shear deformation is much higher when compared to its bending rigidity. Therefore, we neglect the contribution of the transverse shear deformation towards the deformation energy of the fractional-order nonlocal solid in the subsequent analysis.}

By combining the fractional-order nonlinear axial strain in the above equation along with the Euler-Bernoulli displacement field given in Eq.~\eqref{eq: euler_bern_beam}, the axial strain can be recast as:  
\begin{equation}
\label{eq: von_karman_ebt}
    \tilde{\epsilon}_{11}(x_1,x_3) = \tilde{\epsilon}_0(x_1) + x_3\tilde{\kappa}(x_1)
\end{equation}
In the above equation, $\tilde{\epsilon}_0(x_1)$ and $\tilde{\kappa}(x_1)$ denote the fractional-order axial and bending strains, respectively. They are expressed in terms of the mid-plane field variables as:
\begin{subequations}
\begin{equation}
    \tilde{\epsilon}_0(x_1) = D_{x_1}^{\alpha} u_0 (x_1) + \frac{1}{2}\left[D_{x_1}^{\alpha} w_0 (x_1)\right]^2
\end{equation}
\begin{equation}
    \tilde{\kappa}(x_1) = -D_{x_1}^{\alpha} \left[\frac{\mathrm{d}w_0(x_1)}{\mathrm{d}x_1}\right]
\end{equation}
\end{subequations}
Note that we make use of the von-K\'arm\'an definition for geometrically nonlinear strains. Following this approach, large deformations are considered only in the transverse direction. For studies based on this definition of geometrically nonlinear strains, the linear elastic constitutive relations developed in \S 2 can still be employed \cite{ciarlet1980justification}. Therefore, the axial stress in the nonlocal isotropic solid subject to thermoelastic loads may be written from Eq.~(\ref{eq: constt_isot}) to be:
\begin{equation}
    \label{eq: constt_isot_beam}
    \tilde{\sigma}_{11}(x_1,x_3)=E\big(\tilde{\epsilon}_{11}(x_1,x_3)-\alpha_0 \theta(x_1,x_3)\big)
\end{equation}
where $E$ is the Young's modulus of the isotropic solid and $\alpha_0$ is the coefficient of thermal expansion for the isotropic solid, as defined earlier.
Using the above defined fractional-order strain and stress fields, the deformation energy $\mathcal{U}$ of the nonlocal beam is obtained as:
\begin{equation}
    \label{eq: def_energy}
    \mathcal{U} = \frac{1}{2}\int_{\Omega} \tilde{\sigma}_{11}(x_1,x_3) \tilde{\epsilon}_{11}(x_1,x_3) \mathrm{d}V
\end{equation}
where $\Omega$ denotes the volume occupied by the beam. The total potential energy functional of the beam subject to distributed axial ($F_a(x_1)$) and transverse forces ($F_t(x_1)$) acting on the mid-plane, assuming no body forces, is given by:
\begin{equation}
    \label{eq: pot_energy}
    \Pi[\textbf{u}(\bm{x})]=\mathcal{U}-\underbrace{\int_{0}^{L}F_a(x_1) u_0(x_1)\mathrm{d}x_1}_{\text{Work done by axial loads}} - \underbrace{\int_{0}^{L}F_t(x_1) w_0(x_1)\mathrm{d}x_1}_{\text{Work done by transverse loads}}
\end{equation}
We now derive the governing equations and the associated boundary conditions for the thermoelastic response of the nonlocal beam in the strong form by imposing optimality conditions on the above functional. {While we derived the governing equations for the fractional-order solid in \S 2.2.2 following the balance principles for intuitive purposes, here we employ the variational methods due to its ability to treat geometrically non-linear systems.} Before presenting the governing equations, we highlight that the objective of this study is to evaluate the elastic response of a 1D beam when subject to combined thermal and mechanical loads. The thermal load consists of a steady-state temperature distribution applied along the length of the beam on the face at $x_3=\pm h/2$. Thus, the independent variation of temperature field $\delta T$, and thereby $\delta \theta$, is identically zero.

\subsection{Governing equations}
The fractional-order governing equations for the thermoelastic response of geometrically nonlinear and nonlocal beams are obtained using variational principles (i.e. by minimizing the total potential energy given in Eq.~(\ref{eq: pot_energy})). They are given as follows:
\begin{subequations}
\label{eq: all_governing_equations}
\begin{equation}
\label{eq: axial_gde}
        \mathfrak{D}^{\alpha}_{x_1} \mathcal{N}(x_1) + F_a(x_1) = 0 ~~\forall~x_1\in(0,L)
\end{equation}
\begin{equation}
\label{eq: transverse_gde}
        D^1_{x_1}\left[\mathfrak{D}^{\alpha}_{x_1}\mathcal{M}(x_1)\right] + \mathfrak{D}^{\alpha}_{x_1} \left[\mathcal{N}(x_1)D_{x_1}^{\alpha} \left[w_0(x_1)\right]\right] + F_t(x_1) = 0 ~~\forall~x_1\in(0,L)
\end{equation}
\end{subequations}
The corresponding essential and natural boundary conditions are obtained as:
\begin{subequations}
\label{eq: all_BCs}
\begin{equation}
\label{eq: axial_bcs}
    \mathcal{N}(x_1) = 0 ~~\text{or}~~ \delta u_0(x_1) = 0~~\forall~x_1\in\{0,L\}
\end{equation}
\begin{equation}
\label{eq: transverse_moment_bcs}
    \mathcal{M}(x_1)=0~~\text{or}~~\delta\left[ D^1_{x_1} w_0(x_1)\right]=0~~\forall~~x_1\in\{0,L\}
\end{equation}
\begin{equation}
\label{eq: transverse_force_bcs}
    D^1_{x_1} \mathcal{M}(x_1) + \mathcal{N}(x_1) D^1_{x_1} \left[w_0(x_1)\right] =0 ~~\text{or}~~\delta w_0(x_1)=0 ~~\forall~x_1\in\{0,L\}
\end{equation}
\end{subequations}
Note that the detailed steps leading to the above fractional-order nonlinear governing equations extend directly from the geometrically nonlinear analysis of fractional-order beams presented in \cite{sidhardh2020geometrically}, hence they are not provided here. In the above Eqs.~(\ref{eq: all_governing_equations},\ref{eq: all_BCs}), $D^1_{x_1}(\cdot)$ denotes the first integer-order derivative with respect to the axial variable $x_1$. 
Note that the fractional derivative $\mathfrak{D}^{\alpha}_{x_1}(\cdot)$ is defined over the interval $(x_1-l_{B},x_1+l_{A})$ unlike the fractional derivative $D^{\alpha}_{x}(\cdot)$ which is defined over the interval $(x_1-l_{A},x_1+l_{B})$. This change in the terminals of the interval of the Riesz Riemann-Liouville fractional derivative follows from the standard integration by parts technique used to simplify the variational integrals (see \cite{patnaik2019FEM}).
Further, $\mathcal{N}(x_1)$ and $\mathcal{M}(x_1)$ are axial and bending stress resultants defined in the following manner:
\begin{subequations}
\label{eq: stress_resultants}
\begin{equation}
\label{eq: stress_resultant_def_1}
    \mathcal{N}(x_1)=\int_{-b/2}^{b/2}\int_{-h/2}^{h/2}\tilde{\sigma}_{11}(x_1,x_3)~\mathrm{d}x_3~\mathrm{d}x_2
\end{equation}
\begin{equation}
\label{eq: stress_resultant_def_2}
    \mathcal{M}(x_1)=\int_{-b/2}^{b/2}\int_{-h/2}^{h/2}x_3~\tilde{\sigma}_{11}(x_1,x_3)~\mathrm{d}x_3~\mathrm{d}x_2
\end{equation}
\end{subequations}
By using the constitutive relations for a homogeneous isotropic solid given in Eq.~(\ref{eq: constt_isot_beam}) along with the above definitions, the stress resultants are obtained as:
\begin{subequations}
\label{eq: stress_resultants_smp}
\begin{equation}
\label{eq: stress_resultant_smp_1}
    \mathcal{N}(x_1) = A_{11}\tilde{\epsilon}_0(x_1) - N_\theta(x_1)
\end{equation}
\begin{equation}
\label{eq: stress_resultant_smp_2}
    \mathcal{M}(x_1) = -D_{11}\tilde{\kappa}(x_1) - M_\theta(x_1)
\end{equation}
\end{subequations}
where $A_{11}=Ebh$ and $D_{11}=Ebh^3/12$ are the axial and bending stiffness coefficients of the beam, respectively. The thermal resultants $N_\theta(x_1)$ and $M_\theta(x_1)$ for the isotropic beam are given as:
\begin{equation}
    \label{eq: thermal_resultants}
    \left\{N_{\theta}(x_1),M_{\theta}(x_1)\right\}=Eb\alpha_0\int_{-h/2}^{h/2}~\{1,~x_3\}~\theta(x_1,x_3)~\mathrm{d}x_3
\end{equation}
Note that for a general distribution of material properties across the thickness of the beam, additional terms due to the bending-extension coupling would be noted in the Eqs.~(\ref{eq: stress_resultants},\ref{eq: thermal_resultants}). Upon ignoring the nonlinear terms in the governing equations given above, striking similarity may be noted to the governing equations in Eq.~\eqref{eq: momentum_balance_loc} derived using linear momentum balance law.

In the following, we discuss a few characteristics of the thermomechanical governing equations given in Eq.~(\ref{eq: all_governing_equations}). First, observe that the stress resultants given in Eq.~(\ref{eq: stress_resultants_smp}) introduce the thermoelastic variables into the governing equations in Eq.~(\ref{eq: all_governing_equations}). In the absence of thermal loads ($\theta(x_1,x_3)=0$), the constitutive models reduce to the expressions derived in \cite{sidhardh2020geometrically} for the geometrically nonlinear fractional-order nonlocal beams. Owing to the nonlinear nature of the structural response, the axial and transverse displacement fields are coupled unlike what seen in the linear elastic case \cite{patnaik2019FEM}. Second, we emphasize that the fractional-order model for the linear thermoelastic response of a nonlocal beam can be obtained by ignoring the nonlinear terms in the constitutive relations developed above. 
The linear thermoelastic model of the fractional-order nonlocal beam will be discussed further in \S\ref{sec: num_results}. Finally, the classical thermoelastic models are recovered for $\alpha=1$.

\section{Nonlinear fractional finite element model (f-FEM)}
\label{sec: ffem}
Given the nonlinear and integro-differential nature of the governing equations, it is unlikely to obtain closed form solutions for the most general loading and boundary conditions. Therefore, we employ a fractional-order finite element method to obtain the numerical solution of the nonlinear governing equations. The f-FEM developed to solve the thermomechanical fractional-order BVP builds upon the numerical solvers developed for fractional-order models of nonlocal elasticity \cite{patnaik2019FEM,sidhardh2020geometrically}. 
Note that, although the f-FEM is developed and applied for thermoelastic response of fractional-order beams, it can be easily extended to higher dimensional structures like plates and shells.

Analogously to traditional FEM, the f-FEM is formulated starting from a discretized form of the total potential energy functional $\Pi[\textbf{u}(\bm{x})]$ given in Eq.~(\ref{eq: pot_energy}). For this purpose, the 1D domain $\Omega=[0,L]$ of the beam indicated in Fig.~\eqref{fig: schematic_beam} is uniformly discretized into disjoint two noded elements $\Omega^e_i=(x_1^i,x_1^{i+1})$ of length $l_e$ such that $\cup_{i=1}^{N_e} \Omega^e_i=\Omega$, $N_e$ being the total number of discretized elements. It is immediate that $\Omega^e_j\cap\Omega^e_k=\emptyset~\forall~j\neq k$. The unknown field variables $u_0(x_1)$ and $w_0(x_1)$ in Eq.~\eqref{eq: all_governing_equations} can now be evaluated at any point $x_1 \in \Omega^e_i$ by interpolating the corresponding nodal values for $\Omega^e_i$ as:
\begin{equation}
    \label{eq: intpl_disp_u}
    \{u_0(x_1)\}=[\mathcal{L}(x_1)]\{U_e (x_1)\};~~~\{w_0(x_1)\}=[\mathcal{H}(x_1)]\{W_e (x_1)\}
\end{equation}
where $\{U_e (x_1)\}$ and $\{W_e (x_1)\}$ are the axial and transverse displacement degrees of freedom of the two-noded element $\Omega^e_i$. $[\mathcal{L}(x_1)]$ and $[\mathcal{H}(x_1)]$ are the Lagrangian and Hermitian interpolation functions, respectively, chosen to enforce the continuity of the axial and transverse displacement fields for the Euler-Bernoulli beam theory. In terms of the above discussed numerical discretization of the nonlocal domain, the fractional derivatives can be expressed as:
\begin{subequations}
\label{eq: frac_der_full}
\begin{equation}
    D^{\alpha}_{x_1}\left[u_0(x_1)\right]=[\tilde{B}_{u}(x_1)]\{U_g\}; ~ D^{\alpha}_{x_1}\left[w_0(x_1)\right]=[\tilde{B}_{w}(x_1)]\{W_g\}; ~ D^\alpha_{x_1}\left[D^1_{x_1} w_0(x_1)\right]=[\tilde{B}_{\theta}(x_1)]\{W_g\}
\end{equation}
where the matrices $[\tilde{B}_{\square}(x_1)]$ corresponding to nonlocal strain-displacement matrices are given as:
\begin{equation}
\label{eq: frac_der_Bu_final}
    [\tilde{B}_\square(x_1)]=\int_{x_1-l_A}^{x_1+l_B}\mathcal{A}(x_1,s_1,l_A,l_B,\alpha)[B_\square(s_1)][\tilde{\mathcal{C}_\square}(x_1,s_1)]\mathrm{d}s_1
\end{equation}
The kernel $\mathcal{A}(x_1,s_1,l_A,l_B,\alpha)$ in the above equation is: 
\begin{equation}
\label{eq: atten_func}
    \mathcal{A}(x_1,s_1,l_A,l_B,\alpha)=\begin{cases}
    \frac{1}{2}(1-\alpha)l_A^{\alpha-1}{(x_1-s_1)^{-\alpha}}& ~~ s_1\in{(x_1-l_A,x_1)}\\
    \frac{1}{2}(1-\alpha)l_B^{\alpha-1}{(s_1-x_1)^{-\alpha}}& ~~ s_1\in{(x_1,x_1+l_B)}
    \end{cases}
\end{equation}
and the matrices $[B_\square(x_1)]$ $\square \in\{u,w,\theta\}$ are the integer-order strain-displacement matrices given by:
\begin{equation}\label{eq: b_mats}
    [B_u(s_1)]=\frac{\mathrm{d}[\mathcal{L}(s_1)]}{\mathrm{d}s_1};~[B_w(s_1)]=\frac{\mathrm{d}[\mathcal{H}(s_1)]}{\mathrm{d}s_1};~[B_\theta(s_1)]=\frac{\mathrm{d}^2[\mathcal{H}(s_1)]}{\mathrm{d}s_1^2}
\end{equation}
\end{subequations}
Complete details of the steps involved in the numerical evaluation of system matrices following an additional convolution integral in Eq.~\eqref{eq: frac_der_Bu_final} for nonlocal matrices and resolving the singular kernel in Eq.~\eqref{eq: atten_func} for fractional-order derivatives are provided in \cite{patnaik2019FEM,sidhardh2020geometrically}. 

We use the FE approximations of the different fractional-order derivatives to obtain the algebraic governing equations corresponding to the geometrically nonlinear thermoelastic response of the fractional-order beam. In the interest of a more compact notation, the functional dependence of the different physical quantities on the spatial variables will be implied, unless stated to be constant. The first variation of the potential energy function $\Pi[\textbf{u}]$ defined in Eq.~(\ref{eq: pot_energy}) is obtained as:
\begin{equation}
\label{eq: pot_energy_simp}
    \delta \Pi=b\int_{0}^{L} \int_{-h/2}^{h/2} \delta \tilde{\epsilon}_{11}~ \tilde{\sigma}_{11} \mathrm{d}x_3 \mathrm{d}x_1
    -\int_0^L F_t\delta w_0\mathrm{d}x_1-\int_0^L F_a\delta u_0 \mathrm{d}x_1
\end{equation}
By using the strain-displacement relations in Eq.~\eqref{eq: von_karman_ebt} and the stress-resultants in Eq.~\eqref{eq: stress_resultants} we obtain:
\begin{equation}
    \label{eq: pot_energy_simp_2}
    \delta \Pi = \int_{0}^{L} \left\{ \mathcal{N} \left[ D_{x_1}^{\alpha}\left(\delta u_0\right)\right] + \mathcal{N}\left[D_{x_1}^{\alpha} w_0\right] \left[D_{x_1}^{\alpha}\left(\delta w_0\right)\right] - \mathcal{M} \left[D_{x_1}^{\alpha}\left[D^1_{x_1} \left(\delta w_0 \right) \right]\right] - F_a\delta u_0 - F_t \delta w_0 \right\} \mathrm{d}x_1 
\end{equation}
The minimum potential energy principle, $\delta \Pi=0$, is enforced to obtain the algebraic equations of equilibrium. More specifically, by using the numerical approximations developed for the different fractional derivatives (see Eq.~(\ref{eq: frac_der_full})) and then enforcing the minimization of the total potential energy, we obtain the following system of nonlinear algebraic equations in globally assembled vectors of nodal displacement degrees of freedom $\{U_g\}$ and $\{W_g\}$:
\begin{equation}
    \label{eq: algebraic_eqs}
    \begin{bmatrix} [\tilde{K}_{11}] & [\tilde{K}_{12}]\\ [\tilde{K}_{21}] & [\tilde{K}_{22}] \end{bmatrix}
    \begin{Bmatrix} \{U_g\}\\ \{W_g\} \end{Bmatrix}=
    \begin{Bmatrix} \{F_A+F_{A_{\theta}}\}\\ \{F_T+F_{T_{\theta}}\} \end{Bmatrix}
\end{equation}
where the different stiffness matrices are given by:
\begin{subequations}
    \label{eq: stiffness_mats}
\begin{equation}
    \label{eq: stiffness_mats_k11}
    [\tilde{K}_{11}]=\int_{0}^{L}A_{11} [\tilde{B}_{u}(x_1)]^{T} [\tilde{B}_{u}(x_1)] \mathrm{d}x_1
\end{equation}
\begin{equation}
    \label{eq: stiffness_mats_k12}
    [\tilde{K}_{12}]=\frac{1}{2}\int_{0}^{L}\underbrace{A_{11} \left(D_{x_1}^{\alpha}[w_0(x_1)]\right) [\tilde{B}_{u}(x_1)]^{T} [\tilde{B}_{w}(x_1)]}_{\text{Nonlinear matrix}} \mathrm{d}x_1
\end{equation}
\begin{equation}
    \label{eq: stiffness_mats_k21}
    [\tilde{K}_{21}]=\int_{0}^{L}\underbrace{A_{11} \left(D_{x_1}^{\alpha}[w_0(x_1)]\right) [\tilde{B}_{w}(x_1)]^{T} [\tilde{B}_{u}(x_1)]}_{\text{Nonlinear matrix}} \mathrm{d}x_1
\end{equation}
\begin{equation}
    \label{eq: stiffness_mats_k22}
    [\tilde{K}_{22}]=\int_{0}^{L}D_{11} [\tilde{B}_{\theta}(x_1)]^{T} [\tilde{B}_{\theta}(x_1)]~\mathrm{d}x_1+\frac{1}{2}\int_0^L  \underbrace{\left[A_{11}\left(D_{x_1}^{\alpha}[w_0(x_1)]\right)^2\right][\tilde{B}_w(x_1)]^T[\tilde{B}_w(x_1)]}_{\text{Nonlinear matrix}} \mathrm{d}x_1
\end{equation}
\end{subequations}
The axial and transverse force vectors due to the mechanical and thermal loads are given as:
\begin{subequations}
\begin{equation}
    \{F_A\}^{T}=\int_{0}^{L}F_a(x_1) [S_u(x_1)] \mathrm{d}x_1
\end{equation}
\begin{equation}
    \{F_T\}^{T}=\int_{0}^{L}F_t(x_1) [S_w(x_1)] \mathrm{d}x_1
\end{equation}
\begin{equation}
    \{F_{A_\theta}\}^{T}=\int_{0}^{L}N_{\theta}(x_1) [B_u(x_1)] \mathrm{d}x_1
\end{equation}
\begin{equation}
\label{eq: thermal_force_moment}
    \{F_{T_{\theta}}\}^{T}=\int_{0}^{L}N_{\theta}(x_1) \left(D_{x_1}^{\alpha}[w_0(x_1)]\right)~[B_w(x_1)] \mathrm{d}x_1 - \int_{0}^{L}M_{\theta}(x_1) [B_{\theta}(x_1)] \mathrm{d}x_1
\end{equation}
\end{subequations}
{Note that the geometric nonlinearity introduces additional nonlinear thermomechanical coupled terms. This nonlinear behavior is dependent on the thermal properties of the beam as evident from the expressions of the thermal stress resultants given in Eq.~(\ref{eq: thermal_resultants}). In fact, these nonlinear effects are expected to be significant at high temperature. These additional nonlinear thermomechanical terms can be accounted for in two ways: \ul{approach \#1}: the terms are treated as an external nonlinear thermal force; and \ul{approach \#2}: the contribution of these terms is accounted via the stiffness matrix of the system. The equivalence of the results obtained through both these approaches and a comparison of their accuracy and stability is presented in \cite{praveen1998nonlinear}. In this study, we follow the approach \#1 so that the linear analysis of the nonlinear system, for small displacements, becomes straightforward without requiring changes to the stiffness matrix.}

The algebraic equations~(\ref{eq: algebraic_eqs}) are solved for the nodal values of the generalized displacement coordinates for an isotropic beam subject to distributed thermal and mechanical loads. The solution to these equations along with Eq.~(\ref{eq: euler_bern_beam}) gives the displacement field at any point within the beam. The geometric nonlinearity in the system is highlighted by the deformation dependent stiffness terms in Eqs.~\eqref{eq: stiffness_mats_k12}-\eqref{eq: stiffness_mats_k22}. Further, as previously discussed, the additional nonlinear thermomechanical terms are introduced into the model as a nonlinear transverse force as evident from Eq.~(\ref{eq: thermal_force_moment}). Given the nonlinear nature of the FE algebraic equations, a Newton-Raphson (NR) iterative numerical scheme was adopted to obtain the solution of the Eq.~(\ref{eq: algebraic_eqs}). Similar to classical nonlinear models, the NR procedure for the fractional-order nonlinear equations also requires the evaluation of the tangent stiffness matrix. The procedure to evaluate the tangent stiffness matrix as well as the NR scheme can be found in \cite{sidhardh2020geometrically}.

The linear f-FEM for the thermoelastic response of the nonlocal isotropic beam can be obtained from the above model by ignoring the contribution of the nonlinear coupling term, that is $(D_{x_1}^{\alpha}[w_0(x_1)])^2$ in the system matrices as well as in the force vectors. Note that the axial and transverse displacement fields for the linear elastic response due to thermomechanical loads are decoupled. Finally, the f-FEM reduces to a local thermoelastic study of beams when the fractional-order is set to $\alpha=1$.

\section{Numerical results and discussion}
\label{sec: num_results}

\begin{figure}
    \centering
    \includegraphics[width=0.8\textwidth]{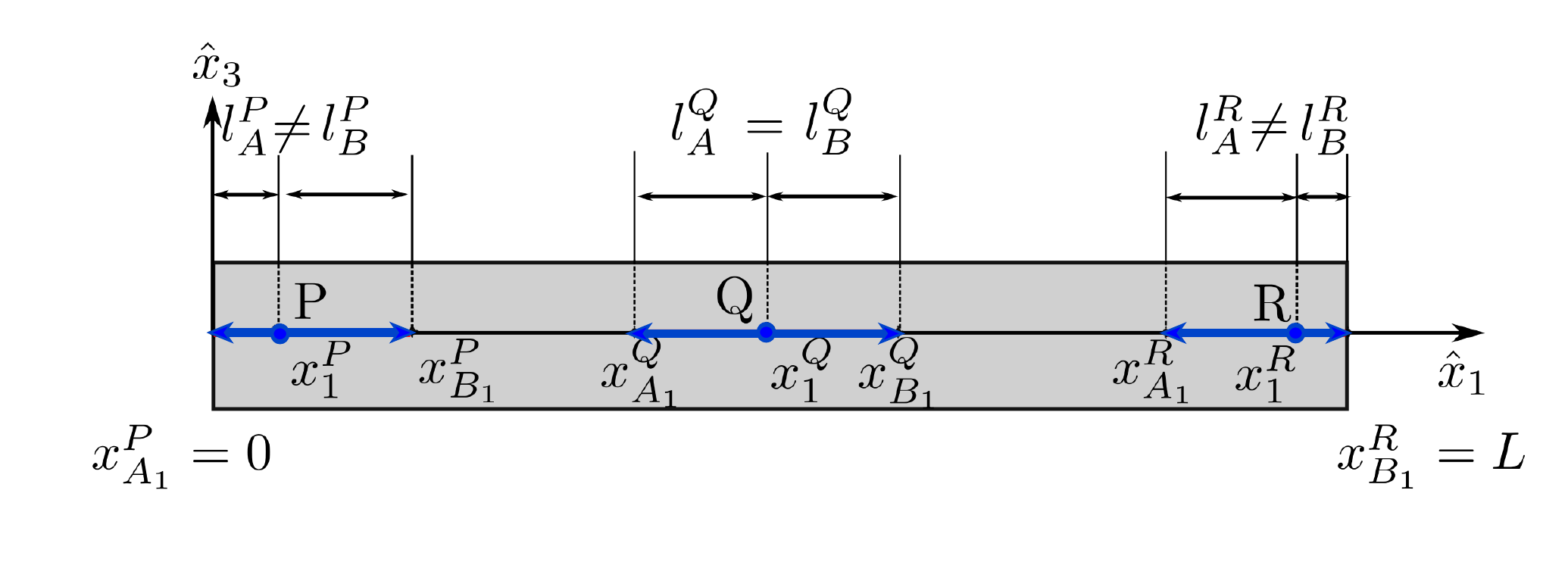}
    \caption{Schematic illustration of the variable nature of the nonlocal length scales for different points along the length of the beam. The length scale on the left side of the point P, $l_A^P$ is truncated such that $l^P_A<l^Q_A$. Similarly, $l^R_B<l^Q_B$.}
    \label{fig: beam_horizon}
\end{figure}

We use the numerical model developed in \S\ref{sec: ffem} to analyze both the linear and the geometrically nonlinear thermoelastic response of fractional-order nonlocal isotropic beams. 
In order to satisfy the underlying assumptions of the Euler-Bernoulli beam theory, the beam is assumed to be slender with an aspect ratio of $L/h=100$. In the following studies, the length of the beam is maintained at $L=1$m and the width of the beam is considered to be unity. The beam is assumed made out of aluminum that is $E=70$GPa and $\alpha_0 = 23\times10^{-6}~{\text{K}}^{-1}$. The constitutive parameters of the fractional-order continuum model, order $\alpha$ and length scales $l_A$ and $l_B$ are provided wherever necessary.
{The length scales $l_A$ and $l_B$ at a point within the domain of the isotropic beam are \textit{considered} equal, that is $l_A=l_B=l_f$.}
However, following the discussion in \cite{patnaik2019FEM}, these length scales are truncated for points close to geometric boundaries of the beam. \SAS{The variable nature of the nonlocal length scales is demonstrated in the schematic given in Fig. \ref{fig: beam_horizon}.}

We analyzed numerically the effect of the fractional parameters $\alpha$ and $l_f$ on the response of the beam subject to different loading and boundary conditions. Both the linear and nonlinear cases are considered. 

Before presenting the results, we make a few remarks concerning the validation and convergence of the f-FEM. In this regard, the f-FEM procedure has already been validated for linear BVPs in \cite{patnaik2019generalized} and for nonlinear BVPs in \cite{sidhardh2020geometrically}. Further, as discussed in \cite{patnaik2019generalized,sidhardh2020geometrically}, the convergence of the f-FEM with finer element discretization is controlled by the “dynamic rate of convergence” defined as: $N_{inf}(=l_f/l_e)$, where $l_e$ is the length of the discretized element. This parameter was shown to be dependent on the fractional model parameters $\alpha$ and $l_f$ that determine the strength of the nonlocal interactions between distant elements \cite{patnaik2019generalized,sidhardh2020geometrically}. Following these convergence studies, the mesh discretization $N_{inf}=10$ was chosen. This choice allows for sufficient number of elements to be included in the horizon of nonlocality at any point, in order to accurately capture the nonlocal interactions \cite{patnaik2019generalized,sidhardh2020geometrically}.\\
 
\begin{figure*}[h!]
    \centering
    \begin{subfigure}[t]{0.5\textwidth}
        \centering
        \includegraphics[width=\textwidth]{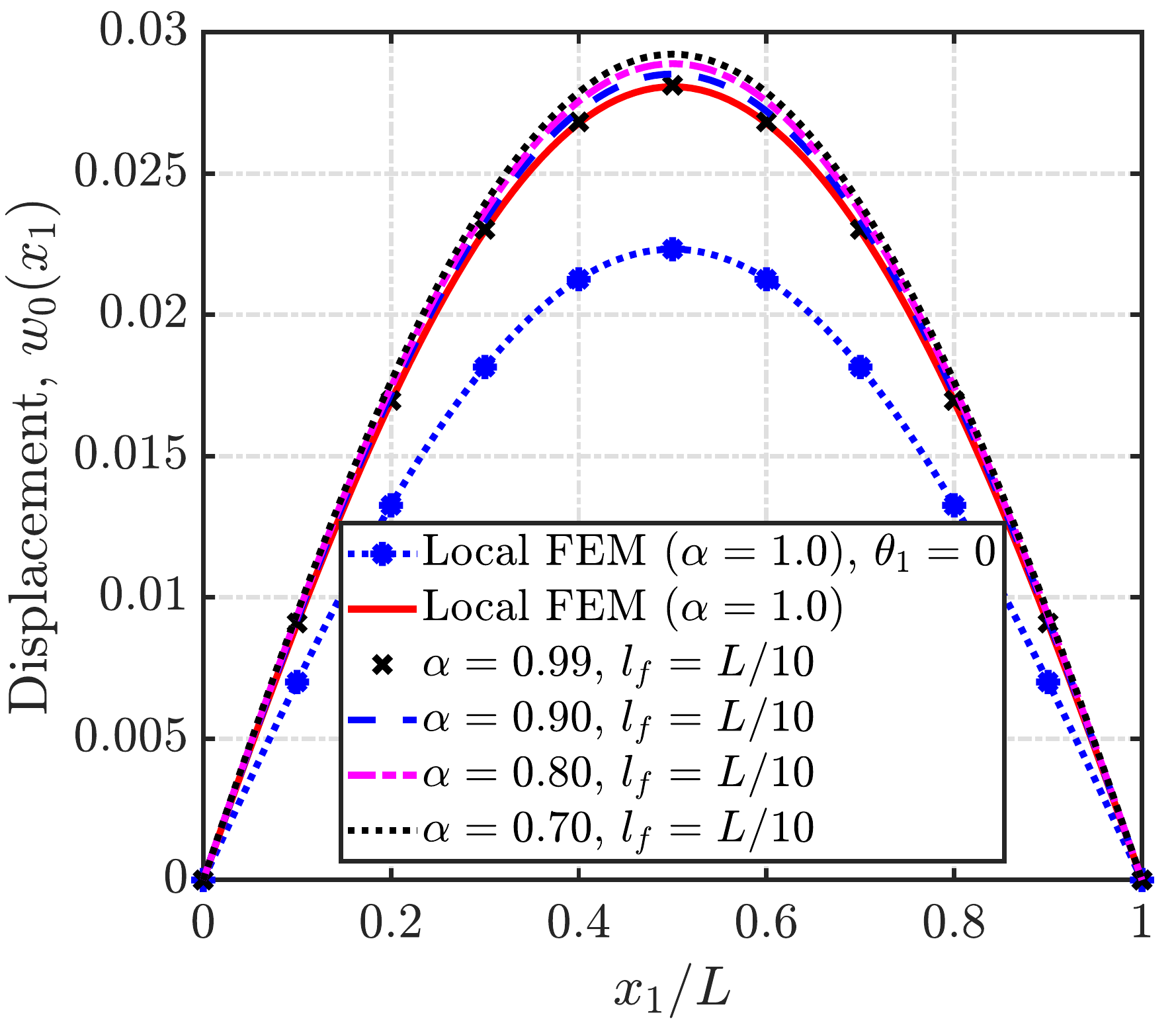}
        \caption{$w_0$ vs $\alpha$.}
        \label{fig: linear_ss_w_alpha}
    \end{subfigure}%
    ~ 
    \begin{subfigure}[t]{0.5\textwidth}
        \centering
        \includegraphics[width=\textwidth]{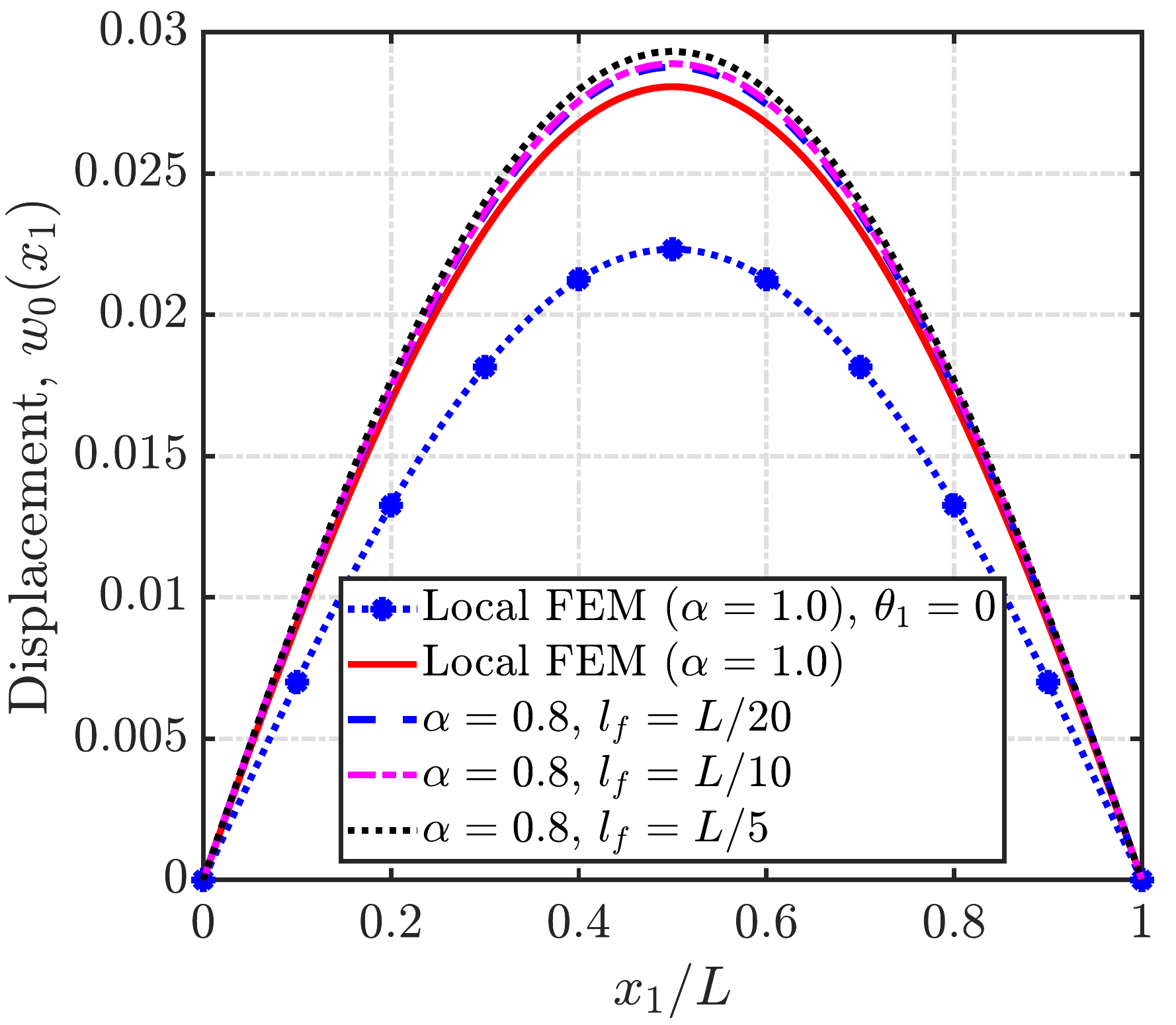}
        \caption{$w_0$ vs $l_f$.}
        \label{fig: linear_ss_w_lf}
    \end{subfigure}
    \caption{Transverse displacement corresponding to the linear response of a simply supported beam for $q_0=10^4$N/m and $\theta_1=10$K. The plot is parameterized for different values of the fractional model parameters.}
    \label{fig: linear_ss_w}
\end{figure*}

\noindent
{\textbf{Linear thermoelastic response}}:
we considered a simply supported beam subject to a uniformly distributed transverse load (UDTL) of magnitude $q_0$ (in N/m) and to the following thermal load
:
\begin{equation}
    \label{eq: thermal_load_linear}
    \theta(x_1,x_3)=\theta_1\left(1+\frac{2x_3}{h}\right)
\end{equation}
{Note that the above distribution is obtained from solving the Fourier's conduction law corresponding to uniform temperatures being applied at the top and bottom suyrfaces ($x_3=\pm h/2$) of the isotropic beam.} Here, the bottom surface of the beam is maintained at the ambient temperature $T_0$. It follows, from Eq.~(\ref{eq: thermal_load_linear}), that the temperature of the top surface is $T_1 ( = T_0 + 2\theta_1)$.  The nonlocal elastic response to these thermomechanical loads for different values of fractional constitutive parameters $\alpha$ and $l_f$ are compared in Fig.~(\ref{fig: linear_ss_w}). The transverse displacement along the length of the simply supported beam was explored for different values of $\alpha$ while maintaining $l_f$ constant (Fig.~(\ref{fig: linear_ss_w_alpha})). Similarly, the transverse displacement was evaluated for different values of $l_f$ while maintaining $\alpha$ constant (Fig.~(\ref{fig: linear_ss_w_lf})). The increase in transverse displacement with the increasing degree of nonlocality, achieved either by reducing $\alpha$ (see Fig.~(\ref{fig: linear_ss_w_alpha})) or by increasing $l_f$ (see Fig.~(\ref{fig: linear_ss_w_lf})), points towards the reduction of the stiffness of the fractional-order beam. We emphasize that the consistent softening of the structure with increasing degree of nonlocality was also observed for beams subject to different loading and boundary conditions. In order to facilitate the analysis of the thermoelastic response of the nonlocal beam, we have provided, as a reference, the transverse displacement of the local beam ($\alpha=1$) for two different cases: (a) absence of thermal load, i.e. $\theta_1=0$; and (b) linear thermal load given in Eq.~(\ref{eq: thermal_load_linear}). From Fig.~(\ref{fig: linear_ss_w}), note  that the nonlocal results converge to the local results for $\alpha$ approaching $1$ and $l_f/L<<1$. 

\begin{figure*}[b!]
    \centering
    \begin{subfigure}[t]{0.5\textwidth}
        \centering
        \includegraphics[width=\textwidth]{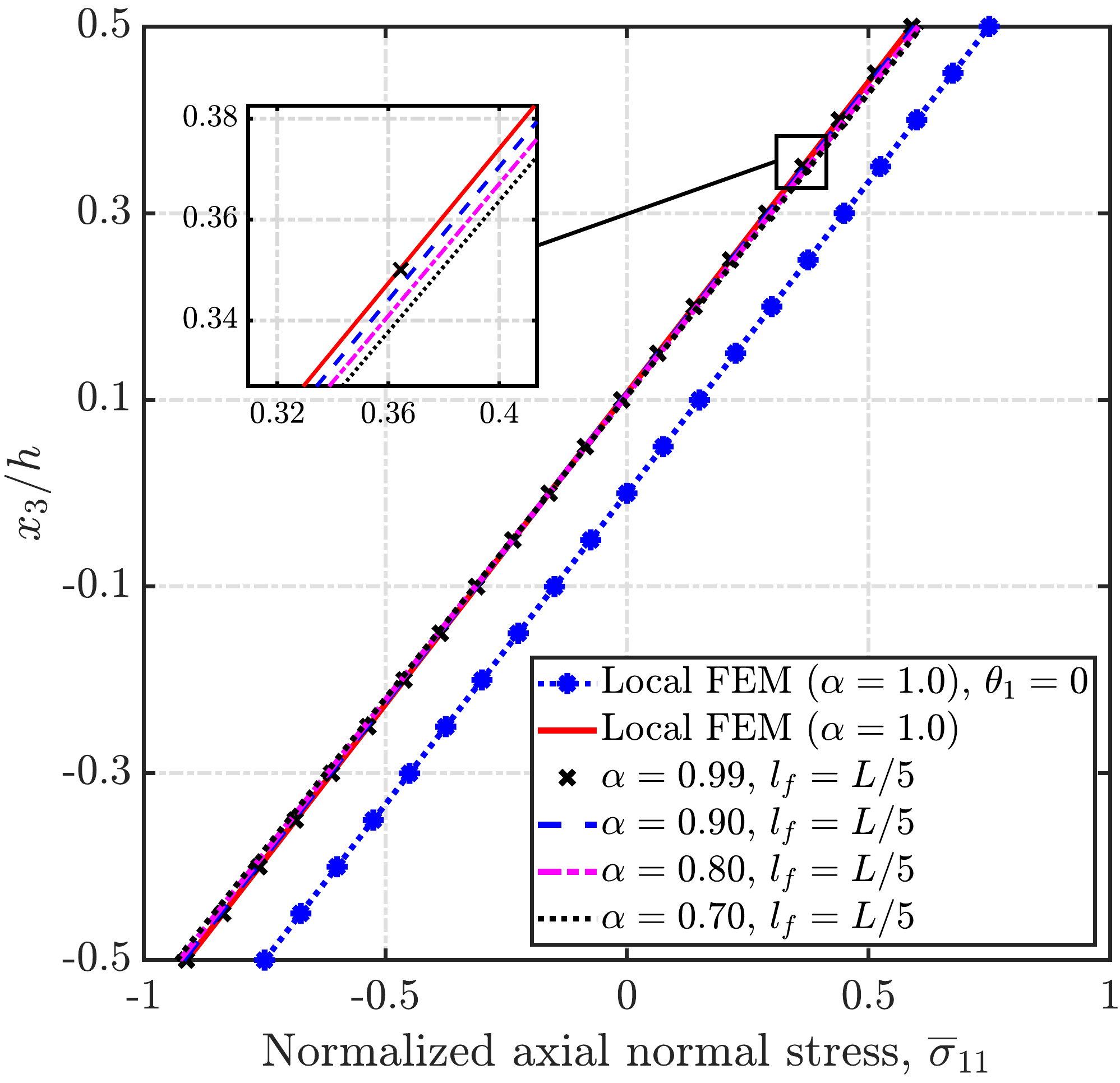}
        \caption{$\overline{\sigma}_{11}$ vs $\alpha$.}
        \label{fig: linear_ss_sig_alpha}
    \end{subfigure}%
    ~ 
    \begin{subfigure}[t]{0.5\textwidth}
        \centering
        \includegraphics[width=\textwidth]{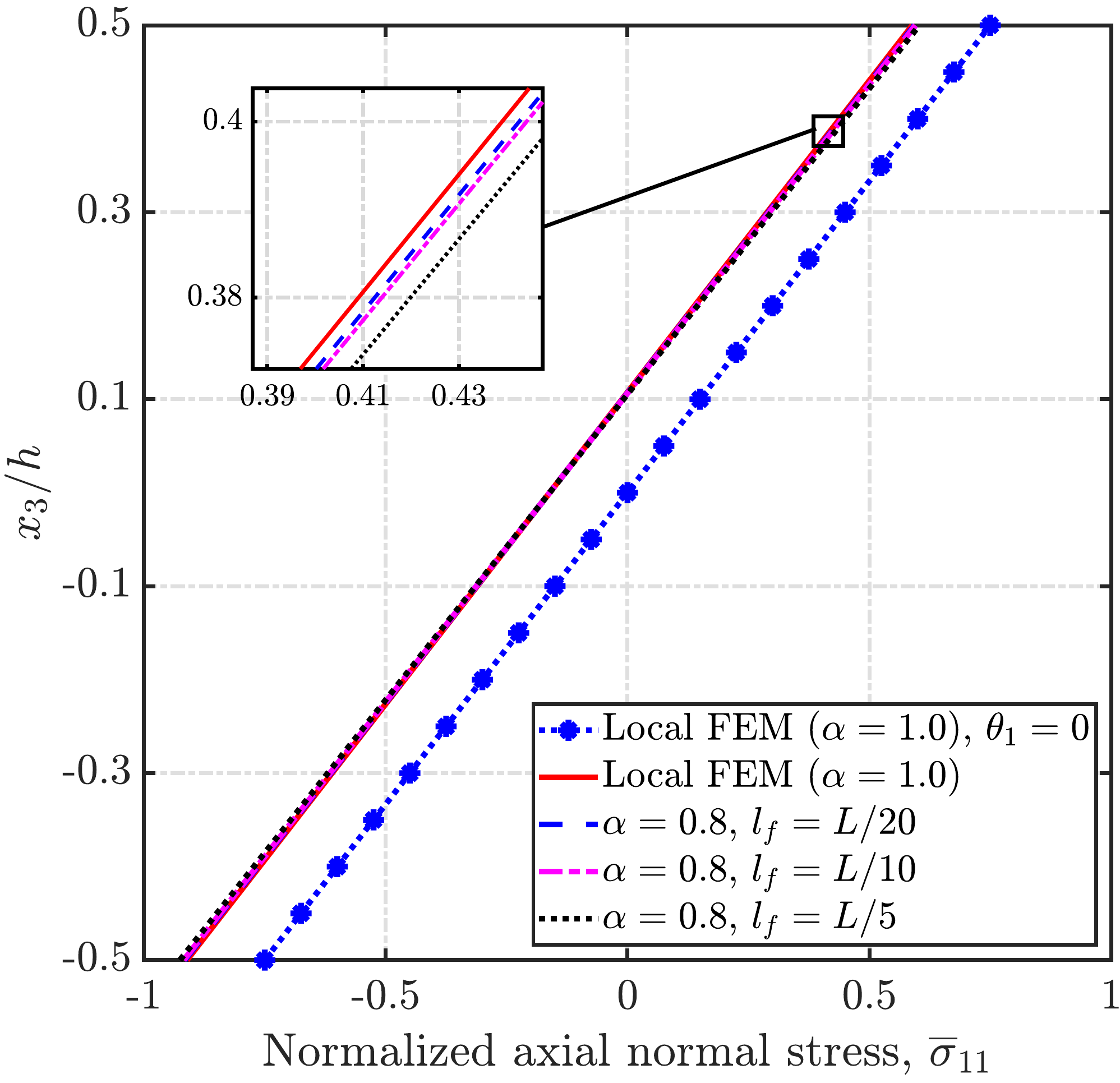}
        \caption{$\overline{\sigma}_{11}$ vs $l_f$.}
        \label{fig: linear_ss_sig_lf}
    \end{subfigure}
    \caption{Normalized values of normal axial stress $\overline{\sigma}_{11}$ across the thickness corresponding to the linear response of a simply supported beam for $q_0=10^4$N/m and $\theta_1=10$K. The plot is parameterized for different values of the fractional model parameters.}
    \label{fig: linear_ss_sig}
\end{figure*}
~
\begin{figure*}[t!]
    \centering
    \begin{subfigure}[t]{0.5\textwidth}
        \centering
        \includegraphics[width=\textwidth]{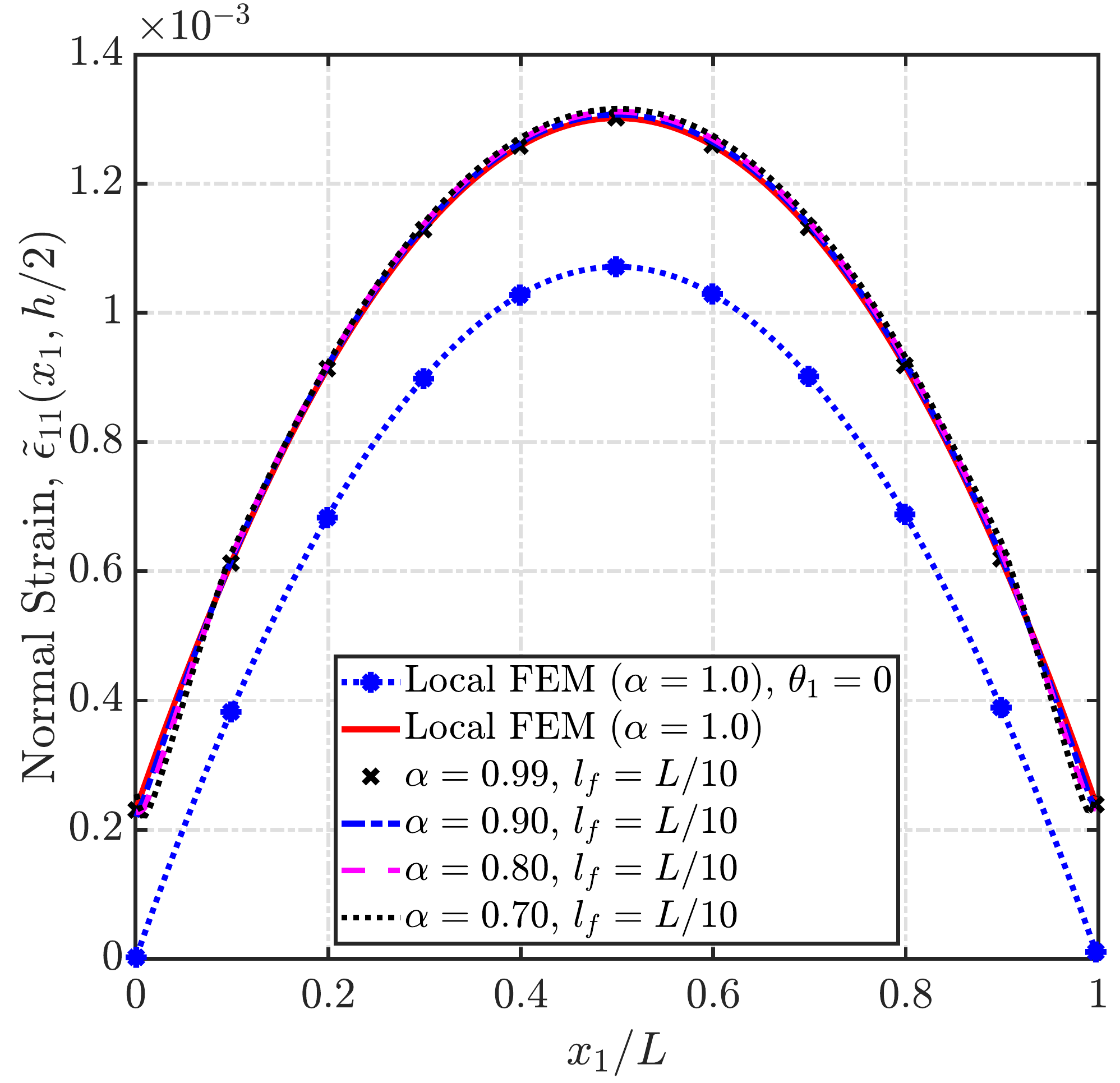}
        \caption{$\epsilon_{11}$ vs $\alpha$.}
        \label{fig: linear_ss_eps_alpha}
    \end{subfigure}%
    ~ 
    \begin{subfigure}[t]{0.5\textwidth}
        \centering
        \includegraphics[width=\textwidth]{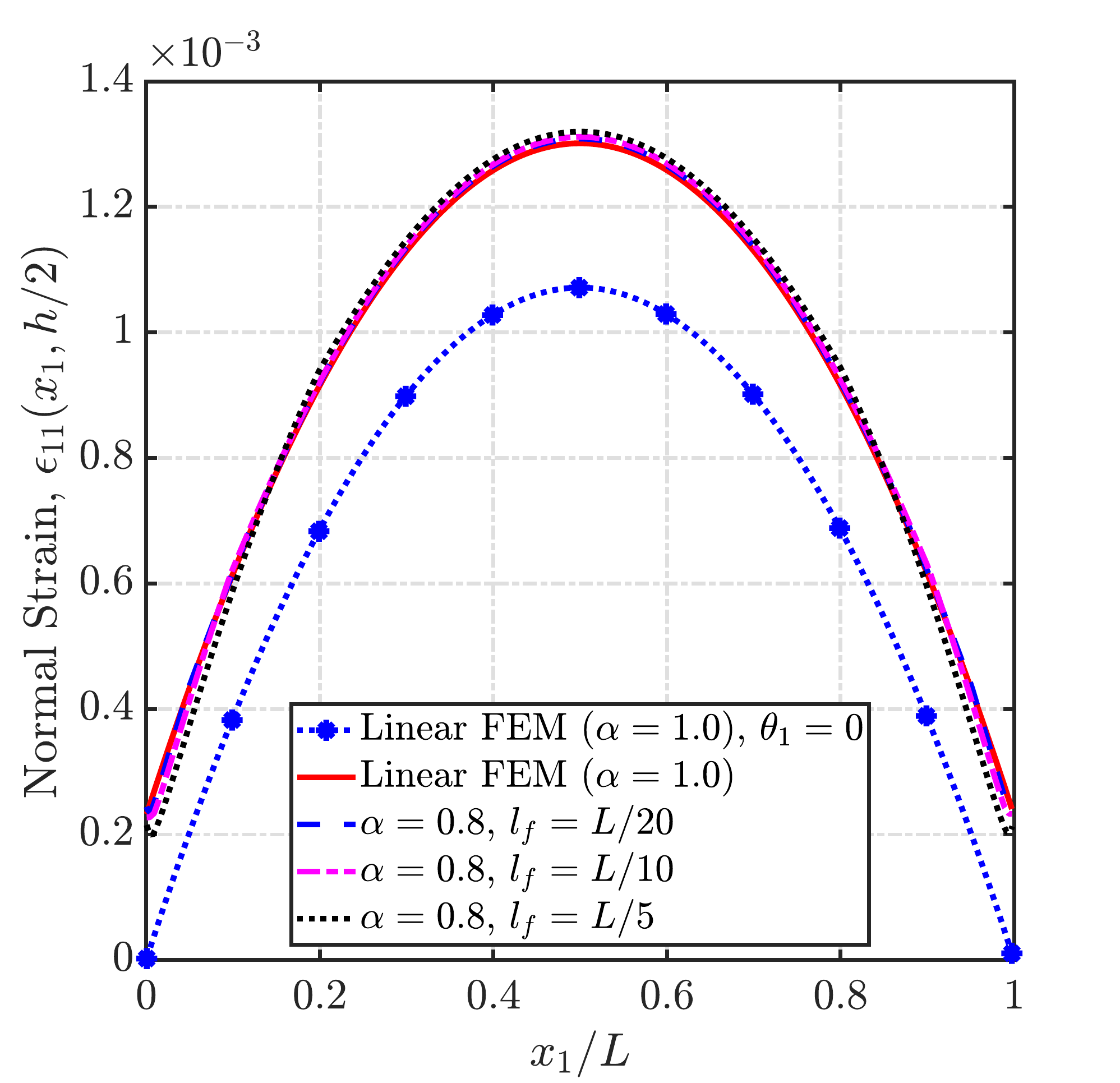}
        \caption{$\epsilon_{11}$ vs $l_f$.}
        \label{fig: linear_ss_eps_lf}
    \end{subfigure}
    \caption{Normal axial strain $\overline{\sigma}_{11}$ corresponding to the linear response of a simply supported beam for $q_0=10^4$N/m and $\theta_1=10$K. The plot is parameterized for different values of the fractional model parameters.}
    \label{fig: linear_ss_eps}
\end{figure*}

\SAS{Extending these studies, the normal axial stress $\sigma_{11}$ and the fractional-order axial normal strain $\epsilon_{11}$ are studied in Figs.~\eqref{fig: linear_ss_sig} and \eqref{fig: linear_ss_eps}. Here, the normal axial stress $\sigma_{11}$ is normalized as follows\cite{khodabakhshi2015unified}:
\begin{equation}
\label{eq: norm_cc}
\overline{\sigma}_{11}(L/2,x_3)=\frac{1}{q_0}\left(\frac{h}{L}\right)^2\tilde{\sigma}_{11}(L/2,x_3)
\end{equation}
In the Fig.~\ref{fig: linear_ss_sig}, (marginally) higher values of $\bar{\sigma}_{11}$ are noted corresponding to an increasing degree of nonlocality. More specifically, an increase in the values of the normalized stress are noted for lower values of the fractional-order $\alpha$ in Fig.~\eqref{fig: linear_ss_sig_alpha}, and higher value of the nonlocal length scale $l_f$ in Fig.~\eqref{fig: linear_ss_sig_lf}. These observations point to a consistent reduction in the stiffness of the fractional-order nonlocal structures, which agrees with the previously noted higher values for the transverse displacement $\bar{w}$ when compared against the response of a local beam (see Fig.~(\ref{fig: linear_ss_w})). Finally, note the constant value of shift in stress profiles at mid-surface caused by the application of the thermal load. Additional results corresponding to the nonlocal axial strain $\tilde{\epsilon}_{11}$ along the length of the beam are presented in Fig.~\ref{fig: linear_ss_eps}.}

\begin{figure*}[t!]
    \centering
    \begin{subfigure}[t]{0.48\textwidth}
        \centering
        \includegraphics[width=\textwidth]{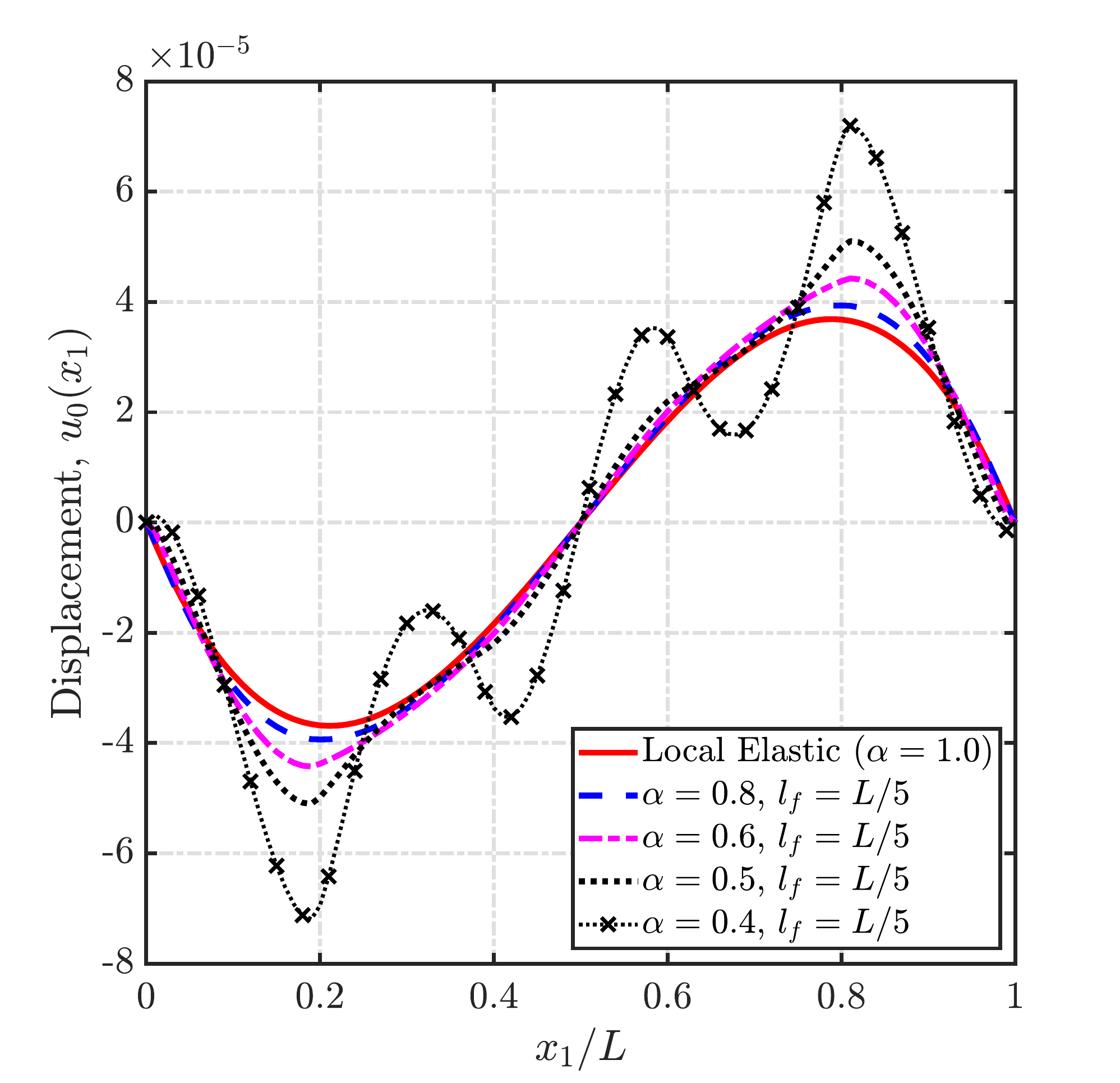}
        \caption{Clamped-Clamped Beam.}
        \label{fig: cc_Axial}
    \end{subfigure}
    ~ 
    \begin{subfigure}[t]{0.48\textwidth}
        \centering
        \includegraphics[width=\textwidth]{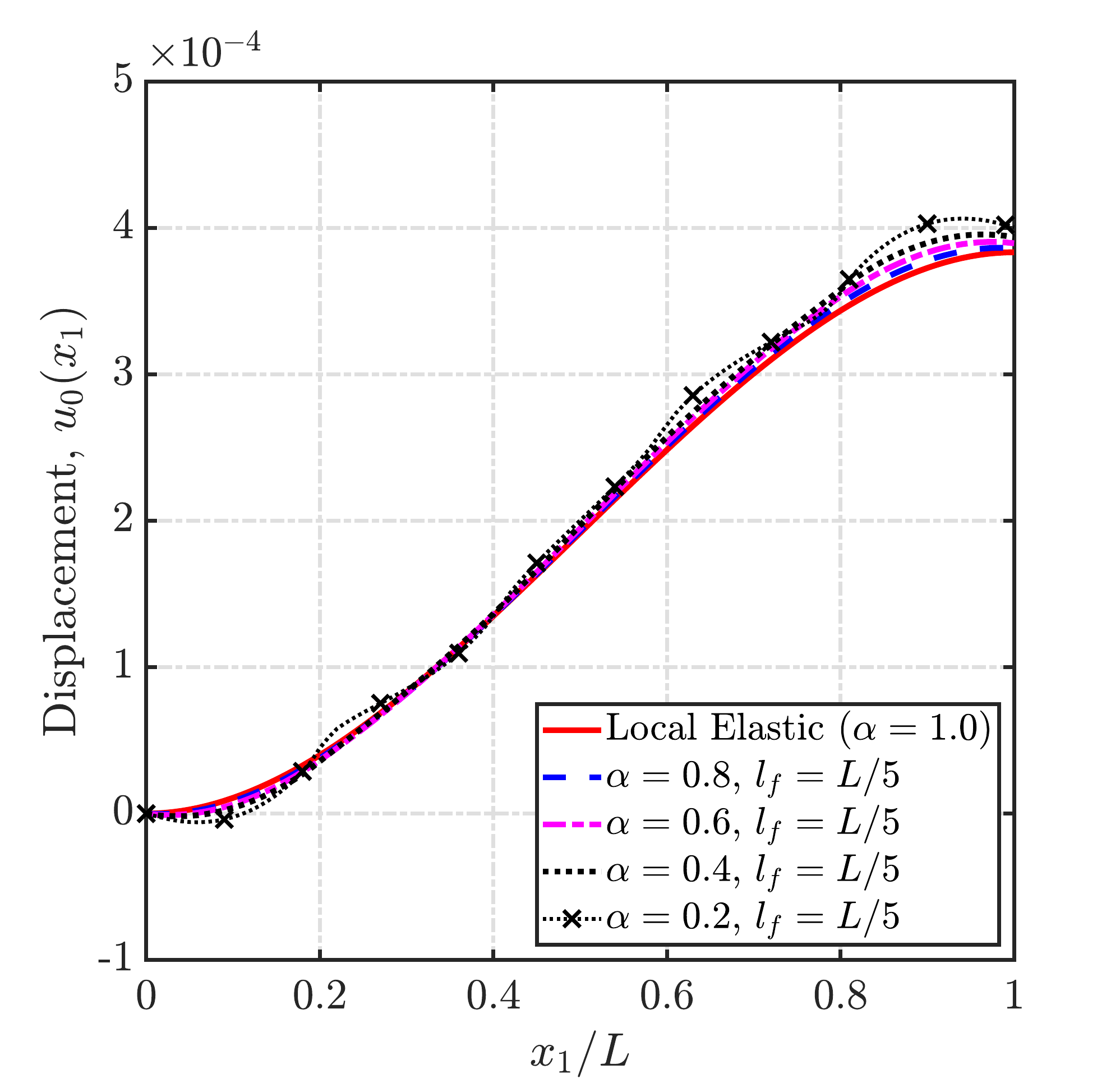}
        \caption{Cantilever Beam.}
        \label{fig: cf_axial}
    \end{subfigure}%
    \caption{Axial displacement corresponding to the linear response of beams subject to (a) clamped-clamped (b) clamped-free (cantilever), boundary conditions for $q_0=0$N/m and a thermal load $\theta_1=100$K. The plot is parameterized for different values of the fractional order with a fixed length scale $l_f=L/5$.}
    \label{fig: comp_theories}
\end{figure*}

In the following, we compare the predictions of the fractional-order approach to nonlocal thermoelasticity with classical integer-order nonlocal theories available in the literature.
More specifically, we compare the predictions of the fractional-order model with the stress-driven integral model presented in \cite{barretta2018stress}. For this purpose, we consider the axial displacement of a doubly-clamped beam subject to the following thermal load in the absence of any mechanical load:
\begin{equation}
\label{eq: thermal_load_length}
    \theta(x_1,x_3)=\theta_1\left(1-\frac{x_1}{L}\right)\frac{x_1}{L}
\end{equation}
The axial displacement of the doubly-clamped beam, obtained via the fractional-order model, is presented in Fig.~(\ref{fig: cc_Axial}) for different values of the fractional-order $\alpha$. As evident from Fig.~(\ref{fig: cc_Axial}), a consistent softening behavior is obtained when considering increasing degree of nonlocality in the fractional model and a constant thermal load. This observation also complements the results obtained in previous studies on fractional-order nonlocal elasticity involving only mechanical loads \cite{patnaik2019FEM}. 
On the contrary, the stress-driven integral approach does not predict either a consistent stiffening or softening response of the same beam for different values of the nonlocal constitutive parameters. In fact, by increasing the degree of nonlocality the beam initially softens and then stiffens, hence leading the authors to label the system as being 'unpredictable'. This result in \cite{barretta2018stress} obtained via the stress-driven integral models for nonlocal thermoelasticity also contrasts with the consistent stiffening predicted by the same model for purely mechanical loads in \cite{romano2017constitutive}. Recall that the fractional-order continuum theories are successful in capturing the softening effects of nonlocal interactions similar to the Eringen's integral theories. The above discussion also highlights an important advantage of the fractional-order approaches to nonlocal elasticity, that is the consistency of the predicted response even in a multi-physics scenario.

{Further, as evident from Figs.~(\ref{fig: linear_ss_w},\ref{fig: comp_theories}), the fractional-order model predicts a reduction in the stiffness of the structure irrespective of the nature of the boundary conditions. This is evident from the consistent softening response shown by the cantilever beams subject to thermal loads given in Eq.~\eqref{eq: thermal_load_length} as depicted in Fig.~(\ref{fig: cf_axial}). This observation of consistent stiffening agrees with the previously studied cases of the doubly clamped and simply supported beams. This is contrary to differential models for nonlocal thermoelasticity, which predict a stiffening response for the cantilever beam and a softening response for other boundary conditions \cite{zenkour2017nonlocal}. As discussed in the introduction, the absence of such paradoxical results in the fractional-order approach to nonlocal elasticity adopted here follows from the positive-definite deformation energy density and self-adjoint nature of the fractional-order nonlocal governing equations \cite{patnaik2019FEM,patnaik2020plates}.} \SAS{Finally, we make an important remark concerning the physically acceptable range for the order $\alpha$. As demonstrated above, the degree of nonlocality increases  with decreasing $\alpha$ leading to a consistent softening of the structure. However, as shown in \cite{sumelka2014fractional,sumelka2015fractional}, results for very low values of $\alpha$ ($\approx0.2$, which indicates a very strong nonlocality in the fractional sense) lead to non-physical solutions. This is also illustrated in the response of fully clamped and cantilever beams subject tot thermal loads in Fig.~\ref{fig: comp_theories}. Hence, there exists a limit on the order of the RC fractional derivative \cite{sumelka2015fractional}. In other terms, it can be concluded that the fractional-calculus based continuum models break down for values of alpha close to the lower integer limit. We emphasize that this breakdown is not a characteristic of the f-FEM technique as the same observation is also noted when using finite difference methods to obtain the numerical solutions (see, for example, \cite{sumelka2014fractional,sumelka2015fractional}). A detailed discussion supported by a possible physical explanation for this loss of consistency in elastic response with reducing fractional order is available in \cite{patnaik2019FEM}. }\\

\begin{figure*}[h!]
    \centering
    \begin{subfigure}[t]{0.5\textwidth}
        \centering
        \includegraphics[width=\textwidth]{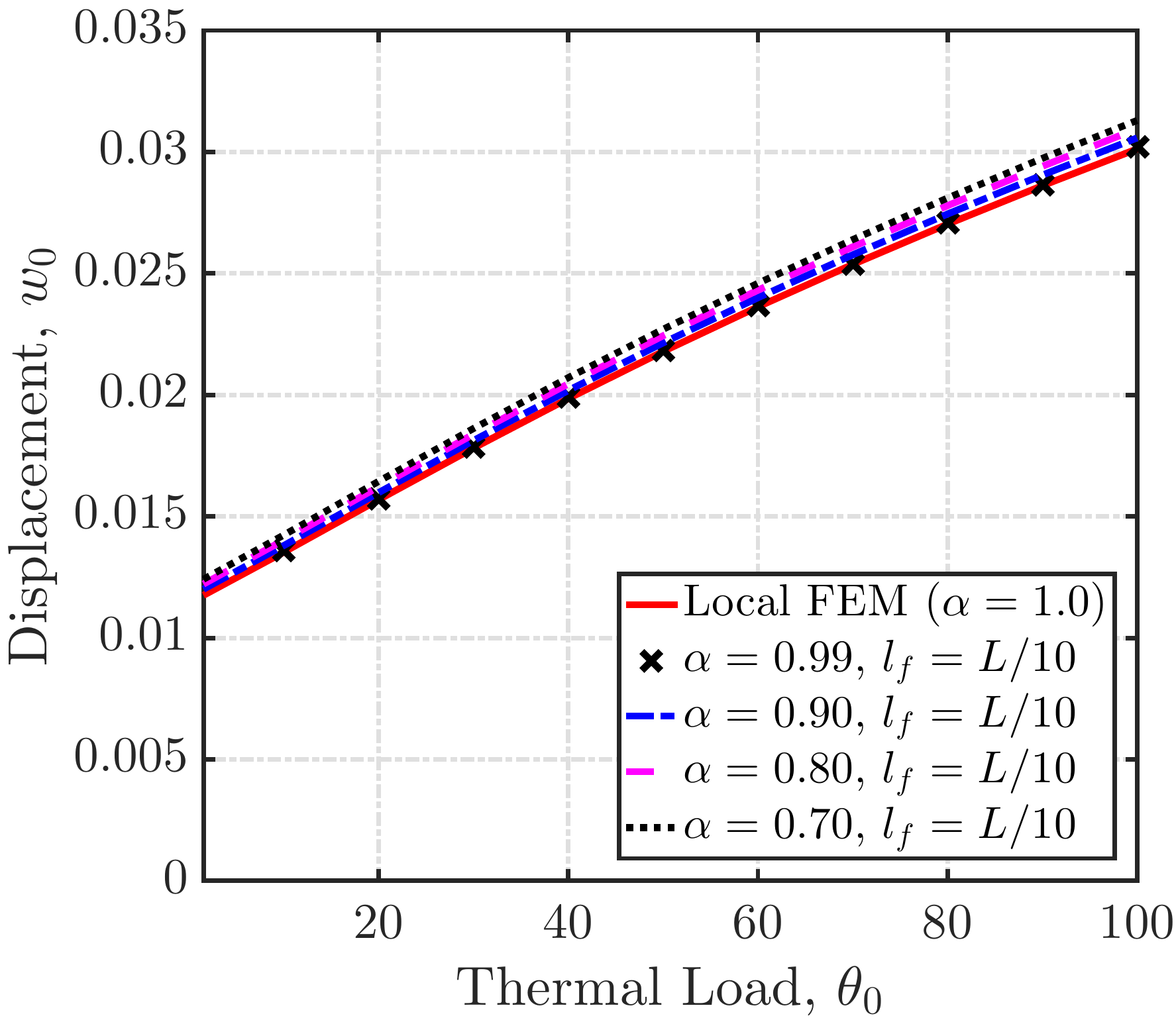}
        \caption{max($w_0$) vs $\alpha$.}
        \label{fig: nonlinear_cc_w_alpha}
    \end{subfigure}%
    ~ 
    \begin{subfigure}[t]{0.5\textwidth}
        \centering
        \includegraphics[width=\textwidth]{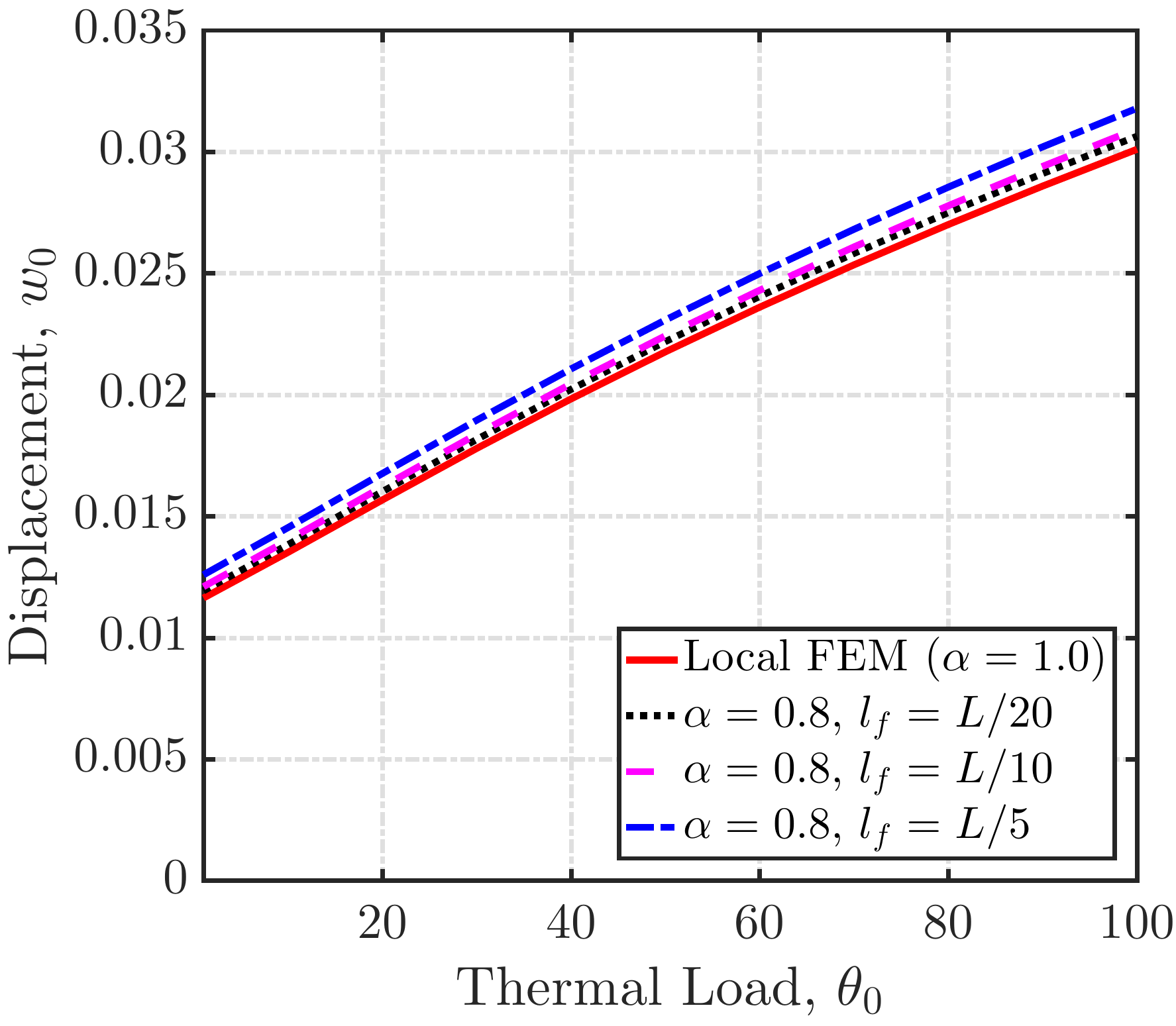}
        \caption{max($w_0$) vs $l_f$.}
        \label{fig: nonlinear_cc_w_lf}
    \end{subfigure}
    \caption{Transverse displacement {at the mid-point} of a clamped-clamped beam subject to $q_0=5\times 10^4$N/m and thermal load $\theta_0$. The curves are parameterized for different values of the fractional-order nonlocal parameters.}
    \label{fig: nonlinear_cc_w}
\end{figure*}

\noindent
\textbf{Nonlinear thermoelastic response}: in this case, the beam was subjected to a UDTL of magnitude $q_0$ (in N/m) and a uniform thermal field $\theta(x_1,x_3)=\theta_0$ (in K) above reference temperature.
First, we considered a beam clamped at both ends and subject to the thermomechanical loads described above. The transverse displacement of the beam for a fixed UDTL and varying magnitude of the thermal load was obtained and compared for different values of $\alpha$ and $l_f$. The magnitude of the UDTL was fixed at $q_0=5\times10^4$N/m and the value of the uniform thermal field $\theta_0$ was varied in order to analyze the effect of the thermal load on the response of the beam. Additionally, in order to analyze the effect of the fractional model parameters on the response of the beam, the transverse displacement of the beam was compared for different values of $\alpha$ and $l_f$. The results of this study are presented in Fig.~(\ref{fig: nonlinear_cc_w}) in terms of the thermal load versus displacement. The displacement values presented in Fig.~(\ref{fig: nonlinear_cc_w}) correspond to the maximum displacement of the mid-plane of the beams, obtained at $x_1=L/2$. The effect of the fractional-order $\alpha$ with $l_f$ being held constant is compared in Fig.~(\ref{fig: nonlinear_cc_w_alpha}). The effect of $l_f$ for fixed $\alpha$ is presented in Fig.~(\ref{fig: nonlinear_cc_w_lf}). The result for the local case ($\alpha=1$) is provided in both cases as a reference. As evident from Fig.~(\ref{fig: nonlinear_cc_w}), the nonlocal beam exhibits consistent softening with increasing thermal loads and increasing degree of nonlocality.
As observed earlier, the thermoelastic response of the nonlocal beam converges to the corresponding local elastic response for $\alpha$ approaching $1$ and $l_f/L \ll 1$.

\begin{figure*}[h!]
    \centering
    \begin{subfigure}[t]{0.5\textwidth}
        \centering
        \includegraphics[width=\textwidth]{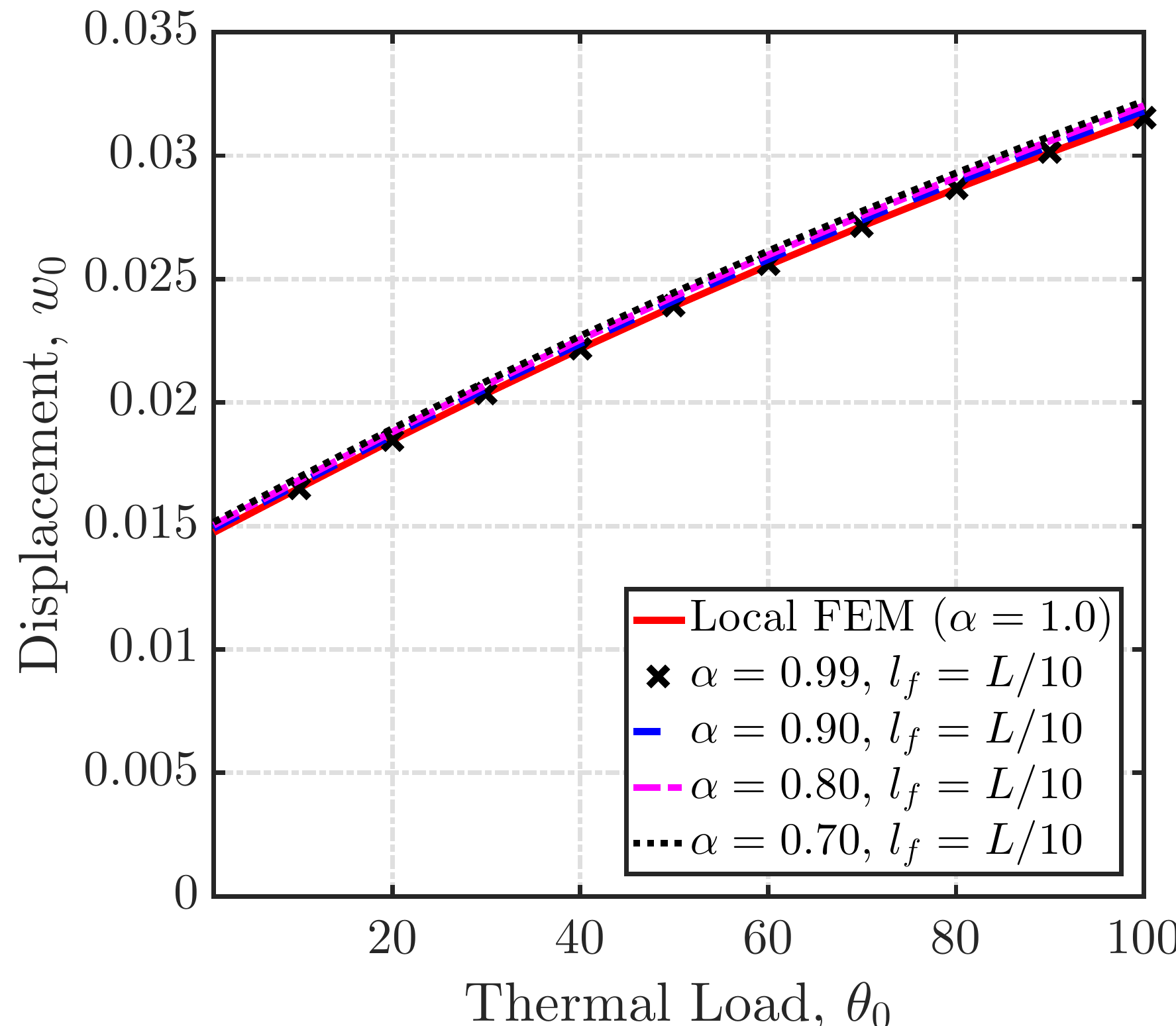}
        \caption{max($w_0$) vs $\alpha$.}
        \label{fig: nonlinear_ss_w_alpha}
    \end{subfigure}%
    ~ 
    \begin{subfigure}[t]{0.5\textwidth}
        \centering
        \includegraphics[width=\textwidth]{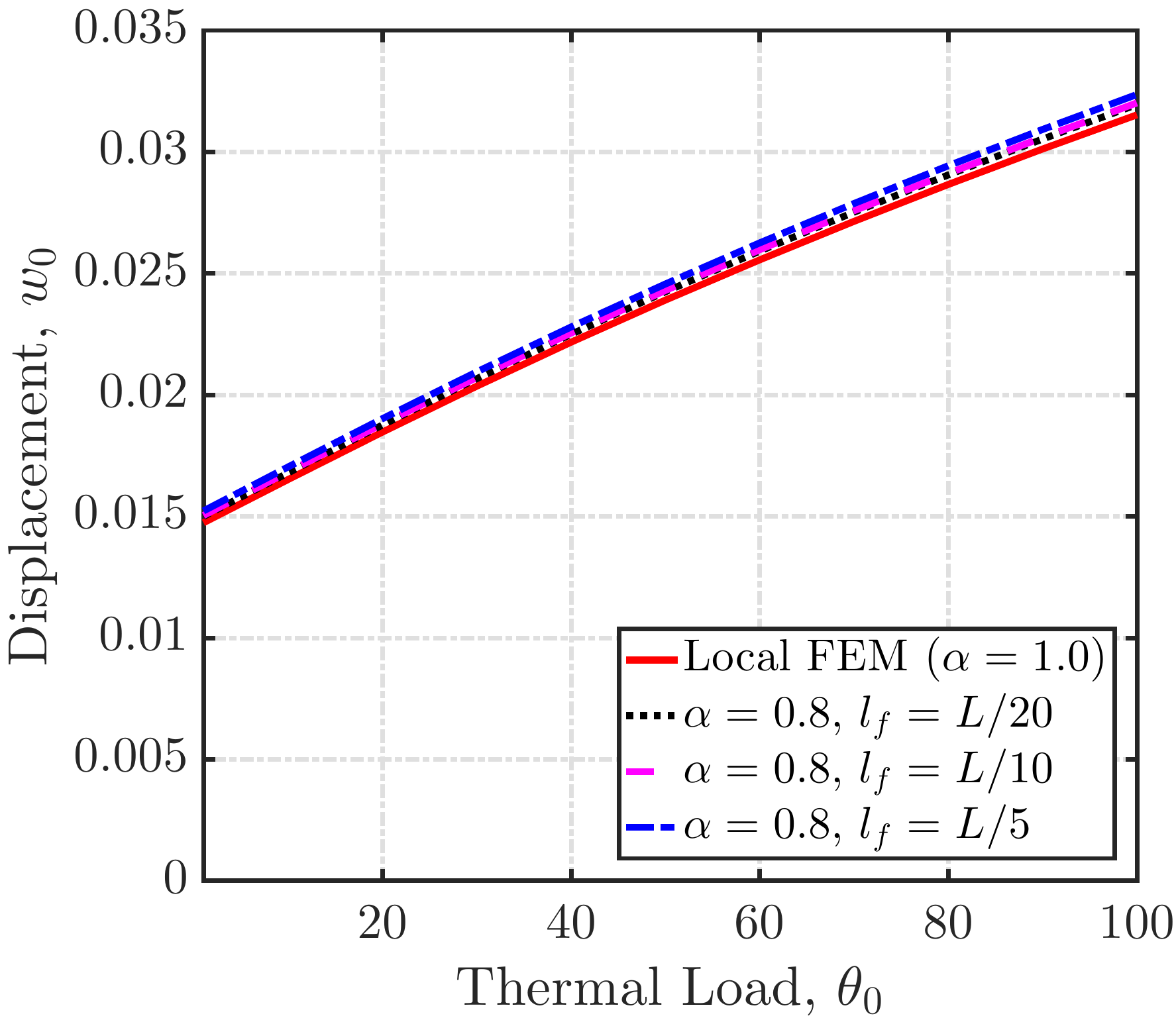}
        \caption{max($w_0$) vs $l_f$.}
        \label{fig: nonlinear_ss_w_lf}
    \end{subfigure}
    \caption{Transverse displacement {at the mid-point} of a pinned-pinned beam subject to $q_0=5\times 10^4$N/m and thermal load $\theta_0$. The curves are parameterized for different values of the fractional-order nonlocal parameters.}
    \label{fig: nonlinear_ss_w}
\end{figure*}

Finally, the study was repeated for a beam pinned at both ends. The effect of the fractional model parameters over the geometrically nonlinear response of the pinned-pinned beam subject to thermomechanical loads is presented in Fig.~(\ref{fig: nonlinear_ss_w}). Observations analogous to those drawn for the clamped-clamped beam can be noted for this case. Remarkably, the fractional-order approach to the modeling of nonlocal elasticity exhibits good consistency across a variety of boundary and loading conditions for both the linear and geometrically nonlinear responses. This behavior differs sharply from the paradoxical results reported in the literature for either gradient or integral based approaches to nonlocal elasticity \cite{challamel2008small,romano2017constitutive}. 

\section{Conclusions}
\label{sec: conclusion}
This study established the thermodynamic framework for fractional-order models of nonlocal thermoelasticity. One of the most significant outcome of this formulation is the ability to rigorously enforce, in a point-wise manner, the thermodynamic balance laws in a nonlocal medium. This result stands in stark contrast with respect to traditional integer-order methods that can satisfy the thermodynamic laws only in a (weak) integral sense. An important consequence of this property of the fractional-order framework is the substantial simplification of the formulation of the free energy density and of the resulting constitutive models of nonlocal thermoelasticity. The thermodynamically-consistent fractional-order continuum theory is well suited to develop accurate models to capture nonlocal interactions, heterogeneity, and scale effects in complex elastic solids operating in a thermomechanical environment. Additionally, the mechanical balance laws are also derived from balance principles. This approach to the development of the governing equations involved a modified Cauchy's lemma for surface tractions in order to include additional forces due to nonlocal effects. The resulting mechanical governing equations are consistent with self-adjoint linear operators admitting unique solutions. This is also an important difference of the present framework in comparison to existing integer- and fractional-order approaches for nonlocal solids available in the literature. 

The efficacy of the fractional-order modeling approach was illustrated by applying the framework to the analysis of the static response of a nonlocal Euler-Bernoulli beam subject to combined thermomechanical loads. Numerical results, obtained using the fractional finite element method, highlighted the extremely robust nature of the fractional models by illustrating the consistency of the predicted nonlocal response across different thermomechanical loads and boundary conditions. A comparison with the existing integer-order nonlocal models was also provided in order to illustrate the advantages offered by the fractional-order theory of nonlocal thermoelasticity. 
\section{Acknowledgements}
The following work was supported by the Defense Advanced Research Project Agency (DARPA) under the grant \#D19AP00052, and the National Science Foundation (NSF) under the grant MOMS \#1761423 and \# DCSD \#1825837. The content and information presented in this manuscript do not necessarily reflect the position or the policy of the government. The material is approved for public release; distribution is unlimited.

\bibliographystyle{unsrt}
\bibliography{thermo_ffem}
\end{document}